\documentstyle[11pt,epsfig]{article}
 
\input xypic

\input{amssymb.sty}

\newtheorem{The}{Theorem}[section]
\newtheorem{Cor}[The]{Corollary}
\newtheorem{Pro}[The]{Proposition}

\newtheorem{Rem}[The]{Remark}
\newtheorem{Lem}[The]{Lemma}

\newtheorem{Def}[The]{Definition}

\def\proof{\vspace{2ex}\noindent{\bf Proof.} }

\def\endproof{\noindent $\diamond$}

\def\dirac{\partial}
\def\d{\partial}

\newcommand{\C}{{\mathbb C}}
\newcommand{\R}{{\mathbb R}}
\newcommand{\Z}{{\mathbb Z}}

\newcommand{\A}{{\cal A}}

\newcommand{\M}{{\cal M}}

\newcommand{\la}{\langle}
\newcommand{\ra}{\rangle}
\newcommand{\ba}{\begin{eqnarray}}
\newcommand{\na}{\end{eqnarray}}
\newcommand{\beq}{\begin{equation}}
\newcommand{\eeq}{\end{equation}} 

\newcommand{\s}{{\mathfrak{s}}}
\newcommand{\spinc}{{\mathrm{Spin}}^c}

\title{Exact triangles in Seiberg-Witten-Floer theory. Part III: proof
of exactness}  
\author{Matilde Marcolli, Bai-Ling Wang}
\date{}

\begin{document}
\maketitle

\tableofcontents

\section{Introduction}

In the previous parts of this work \cite{CMW}, \cite{MW2}, we
considered a homology 3-sphere $Y$ 
and an embedded knot $K$. We considered manifolds $Y_1$ and $Y_0$ obtained by
1-surgery and 0-surgery on $K$, respectively. We proved that the
oriented moduli space of solutions of the Seiberg--Witten equations on
$Y$ can be described as the union 
\begin{equation}
\label{mod:decomp}
{\cal M}_{Y,\mu} \cong {\cal M}_{Y_1} \cup \bigcup_{\s_{k}}  {\cal
M}_{Y_0}(\s_{k}), 
\end{equation}
where $\mu$ is a suitable perturbation that simulates the effect of
surgery, and the $\s_{k}$ are the $\spinc$ structures on $Y_0$ that
reduce to the trivial structure when restricted on the tubular
neighbourhood of the knot $\nu(K)$ and on the knot complement $V$.
The irreducible components are similarly related
$$ {\cal M}^*_{Y,\mu} \cong {\cal M}^*_{Y_1} \cup \bigcup_{\s_{k}}  {\cal
M}_{Y_0}(\s_{k}), $$
as shown in \cite{CMW}.

We also proved that there is a way of assigning compatible choices of
the relative grading, so that (\ref{mod:decomp}) gives an exact
sequence {\em of abelian groups}
\begin{equation}
\label{seq:groups}
0\to C_q (Y_1)\stackrel{j_q}{\to} C_q(Y,\mu) \stackrel{\pi_q}{\to}
\bigoplus_{\s_{k}} C_{(q)} (Y_0,\s_{k}) \to 0.
\end{equation}

Recall that the grading of $C_{(q)} (Y_0,\s_{k})$ denotes a lifting
to a $\Z$-graded   
complex of the $\Z_{2i_{k}}$-graded complex. This lifting is
determined by a compatible 
choice of grading on (\ref{mod:decomp}). This $\Z$-graded Floer
complex is analyzed in 
detail in \cite{MW4}. The maps $j_*$ and $\pi_*$ are
induced by the inclusion
$$ j: {\cal M}_{Y_1} \hookrightarrow {\cal M}_{Y,\mu}, $$
and the projection on the quotient, generated by
the elements of $$\cup_{\s_{k}}  {\cal M}_{Y_0}(\s_{k}).$$

The analysis of the splitting and gluing of the moduli spaces of flow 
lines in \cite{MW2} shows that, in general, one should not expect
(\ref{seq:groups}) to be an exact sequence {\em of chain complexes}:
the maps $j_*$  and $\pi_*$ need not commute with
the boundary operators of the Floer complexes. Thus, on the algebraic
point of view, the existence of the decomposition (\ref{mod:decomp})
simply signifies that the existence of an exact sequence is possible:
in fact, the decomposition (\ref{seq:groups}) shows that the ranks are
compatible with the existence of the desired exact sequence.
The maps that provide the exact sequence are derived from the surgery
cobordisms.

In this paper we introduce maps $w^1_*$ and $w^0_*$, induced 
by the surgery cobordisms $W_1$ and $W_0$ connecting $Y_1$ and $Y$, and $Y$
and $Y_0$, respectively. The maps $w^1_*$ and $w^0_*$ are defined by a 
suitable choice of the $\spinc$ structures $\s_\ell$ and $\s_k$
on the 4-manifolds $W_1$ and $W_0$, that restrict to the assigned
$\spinc$ structure at the two ends of the cobordism. 

We show that $w^1_*$ and $w^0_*$ are chain homomorphisms,
\begin{equation}
C_q (Y_1)\stackrel{w^1_q}{\to} C_q(Y,\mu)
\label{cob:map:1}
\end{equation}
and
\begin{equation} 
C_q(Y,\mu) \stackrel{w^0_q}{\to} \bigoplus_{\s_{k}} C_{(q)} (Y_0,\s_{k}). 
\label{cob:map:2}
\end{equation}

We then show that the map $w_*^0$ is surjective and the map $w_*^1$ is
injective. The main technique is similar to the technique 
developed in \cite{MW2} in order to study the behavior of flow lines
under the splitting $Y=V\cup_{T^2} \nu(K)$. Here we consider punctured
surgery cobordisms $W_1 \backslash \{ x_1 \}$ and $W_0\backslash \{
x_0 \}$ and we stretch product regions $T^2 \times [-r,r] \times \R$
inside these punctured cobordisms, thus also stretching product regions
$T^2\times [-r,r]\times [0,\infty)$ in the $S^3\times [0,\infty)$ end
near the puncture $x_i$. We have geometric limits and a gluing theorem 
as in \cite{MW2} and we can show that, when the parameter $\epsilon$
in the surgery perturbation is small enough, the solutions on $W_1$
and $W_0$ have the following behavior.
If we consider solutions on $W_1$ with asymptotic values
$a_1$ in $\M_{Y_1}$ and $j(a_1)$ in $j( \M_{Y_1})\subset \M_{Y,\mu}$,
and solutions on $W_0$ with asymptotic values $a \in
\M_{Y,\mu}\backslash j( \M_{Y_1})$ and $\pi(a)$ in $\cup_k
\M_{Y_0}(\s_k)$, then the corresponding components of the maps $w^1_*$ and
$w^0_*$ agree with those of the maps $j_*$ and $\pi_*$ in
(\ref{seq:groups}). In other words, we prove the following relations:
\ba \la j(a_1') , w^1_*(a_1) \ra =\delta_{a_1' , a_1}, \label{delta1} \na
for all $a_1'$ and $a_1$ in $\M_{Y_1}$, and
\ba \la a_0 ,w^0_*(a) \ra = \delta_{a_0,\pi(a)}, \label{delta2} \na
for all $a\in \M_{Y,\mu}\backslash j( \M_{Y_1})$ and for all $a_0 \in
\cup_k \M_{Y_0}(\s_k)$.

We then prove the relation $w^0_* \circ w^1_*=0$. Again, we follow the 
same technique. We compare the geometric limits of zero-dimensional
moduli spaces 
$$ \M^{W_1}_\ell (a_1,a) $$
and
$$ \M^{W_0}_k(j(a_1),\pi(a)), $$
for $a\in \M_{Y,\mu}\backslash j(\M_{Y_1})$, regarded as subsets of
the moduli space
$$ \M^{W}_{\ell,k}(a_1, \pi(a)) $$
on the composite cobordism $W=W_1 \cup_Y W_0$.
We obtain an orientation reversing diffeomorphism
$$ \M^{W_1}_\ell (a_1,a) \cong \M^{W_0}_k(j(a_1),\pi(a)), $$
which proves the relation $w^0_* \circ w^1_*=0$. 
The main technique consists of identifying the moduli spaces of
solutions on the cobordisms with certain pre-gluing data obtained out
of the explicit description of the geometric limits developed in Part
II \cite{MW2}. The set of pre-gluing data and the resulting moduli
space can be identified up to a diffeomorphism given by the gluing map.
The set of pre-gluing data consists of moduli spaces of finite energy
monopoles on $V\times \R$ together with holomorphic triangles in a
covering of the character variety $\chi(T^2)$ of flat
$U(1)$-connections on $T^2$, with boundary along arcs of Lagrangians
$\ell$, $\ell_1$, $\ell_\mu^*$ determined by the asymptotic values 
$\partial_\infty \M_V^*$, and by the flat connections on $\nu(K)$,
with or without surgery perturbation. In order to use dimensional
arguments, we compare the formulae for the Maslov index in the
splitting of the spectral flow, with the formulae for the dimension of
the moduli space of such holomorphic triangles.

Then we can complete the proof of the exactness of the sequence 
\begin{equation}
0\to C_q (Y_1)\stackrel{w^1_q}{\to} C_q(Y,\mu) \stackrel{w^0_q}{\to}
\bigoplus_{\s_{k}} C_{(q)} (Y_0,\s_{k}) \to 0.
\label{seq:exact}
\end{equation}
It is enough to show that (\ref{delta1}) and
(\ref{delta2}), together with the inclusion $Im(w^1_*)\subset
Ker(w^0_*)$ determine enough relations among the coefficients of the
chain maps that force the reverse inclusion $Ker(w^0_*)\subset
Im(w^1_*)$ to hold as well.

In the last Section we analyze the connecting homomorphism in the
long exact sequence
$$ \cdots \to HF_q(Y_1)\to HF_q(Y,\mu)\to \oplus_k
HF_{(q)}(Y_0,\s_k)\stackrel{\Delta_q}{\to} HF_{q-1}(Y_1) \to \cdots.$$
We first show that the coefficients of $\Delta$ are given by the
component of the boundary on
$Y$ that counts flow lines connecting critical points of
$\M_{Y,\mu}\backslash j( \M_{Y_1})$ to critical points of
$j( \M_{Y_1})$, of relative index one. We then proceed to identify
this counting with the counting of zero dimensional moduli 
spaces on another cobordism $\bar W_2$ connecting $Y_0$ and $Y_1$,
satisfying the relation
$$ \bar W = \bar W_2 \# \C P^2, $$
where $\bar W$ is the composite cobordism $\bar W=\bar W_0 \cup_Y \bar
W_1$. Thus, we obtain the result that the
exact triangle for Seiberg-Witten Floer homology  
is a surgery triangle, that is, the connecting homomorphism in the
exact sequence is determined by a chain map $\bar w^2_*$ induced by the
surgery cobordism $\bar W_2$, and the resulting diagram
$$ C_*(Y_1) \stackrel{w_*^1}{\to} C_*(Y,\mu) \stackrel{w^0_*}{\to}
\oplus_k C_{(*)}(Y_0,\s_k) \stackrel{\bar w^2_*}{\to} C_*(Y_1)[-1] $$
is a distinguished triangle.

\noindent{\bf Acknowledgments.} The first author is partially supported
by NSF grant DMS-9802480. The second author is partially supported by
ARC Fellowship. We thank the Max--Planck--Institut f\"ur
Mathematik, where a large part of the work was done.

\section{The cobordisms}

We describe briefly the topology of the cobordisms and then introduce
the appropriate perturbed Seiberg--Witten equations, and the
corresponding moduli spaces.

The cobordism $W_1$ is obtained by removing from the trivial cobordism 
$Y_1\times I$ an $S^1\times D\cong \nu(K)\times \{ 1\}$, where $D$ is
a disk, and $\nu(K)$ is the tubular neighbourhood of the knot in
$Y_1$, and then attaching a 2-handle with framing $-1$.
We denote by $D_1$ the core disk of the 2-handle in $W_1$. 
Similarly, the cobordism $W_0$ is obtained by removing
from the trivial cobordism $Y_0\times I$ an $S^1\times D\cong
\nu(K)\times \{ 0\}$ and attaching a 2-handle with framing zero.
We denote by $D_0$ the core disk of the 2-handle in $W_0$.
Attaching the two-handle has the effect of modifying the
boundary component $Y_1\times \{ 1\}$ in the trivial cobordism to the
boundary component $Y\times \{ 1 \}$ in the non-trivial cobordism
$W_1$, or, respectively,  the
boundary component $Y_0\times \{ 0 \}$ in the trivial cobordism 
to the boundary component $Y\times \{ 0 \}$ in $W_0$. 

\begin{Lem}
The cobordisms $W_1$ and $W_0$ have the following topology:
$$ \begin{array}{cccc} b_1(W_1) =0 & b_2(W_1)=1 & b_2^+(W_1)=0 &
b_2^-(W_1)=1 \\
b_1(W_0) =0 & b_2(W_0)=1 & b_2^+(W_0)=0 &
b_2^-(W_0)=0. \end{array} $$
The composite cobordism $W=W_0 \cup_{Y} W_1$, connecting $Y_1$ and
$Y_0$, can be written as a blow up $W= W_2 \# \overline{\C P^2}$, where
$W_2$ satisfies 
$$ \begin{array}{cccc}  b_1(W_2)=0 & b_2(W_2)=1 & b_2^+(W_2)=0 &
b_2^-(W_2)=0. \end{array} $$
\end{Lem}

\proof
Let $\sigma$ be a Seifert surface for the knot $K$ in the homology
sphere $Y$, $\partial \sigma =K$. In the cobordism $W_1$ consider the
surface $\Sigma_1$ obtained by attaching the Seifert surface and the
core disk $D_1$ along the knot $K$,
$$ \Sigma_1 =\sigma\cup_K D_1. $$
The homology $H_2(W_1, \Z)$ is generated 
by the class $[\Sigma_1 ]$ with self intersection $-1$. Similarly,
consider the surface $\Sigma_0$ in $W_0$, obtained by attaching along
$K$ the 
Seifert surface and the core disk $D_0$,
$$ \Sigma_0= \sigma\cup_K D_0. $$
The homology $H_2(W_0, \Z)$ is generated by the class $[\Sigma_0 ]$
with self 
intersection zero. 
The surface 
$$ \Sigma = D_1 \cup_{K} \overline{D_0} $$ 
in the composite cobordism $W$ has self intersection $-1$. The
homology $H_2(W, \Z)$ is generated by $[\Sigma]=[\Sigma_1]
-[\Sigma_0]$ and $[\Sigma_0]$. The class $[\Sigma]$ represents the
exceptional divisor $E$ in the blowup.
Thus, the blown down cobordism $W_2$ has homology $H_2(W_2, \Z)$
generated by a class, which we still write  $[\Sigma_0]$, 
with self intersection zero.

\endproof

Notice that the surface $\Sigma_0$ is homologous in $W_0$ to the
generator $\sigma \cup_K D$ of $H_2(Y_0, \Z)$. The following simple
calculation is useful in classifying the possible
$\spinc$ structures on the cobordisms.

\begin{Lem}
We have $H_2(W_0, Y \cup Y_0,\Z)=\Z$ generated by a surface $\tilde
\Sigma_0$ in $W_0$ with $\partial \tilde\Sigma_0=\gamma$, with $[\gamma]$
the generator of $H_1(Y_0,\Z)$.
The class $[\Sigma_0]$ introduced above, which
generates $H_2(W_0,\Z)$ is trivial in $H_2(W_0, Y \cup Y_0)$.

Moreover, we have $H_2(W_1, Y_1 \cup Y,\Z)=\Z$ generated by the image of
the class $[\Sigma_1]$ which generates $H_2(W_1,\Z)$, and $H_2(W, Y_1
\cup Y_0,\Z)=\Z\oplus \Z$, generated by the class $E$ of the exceptional
divisor of the blow-up and by $[\tilde \Sigma_0]$. Again, the class
$[\Sigma_0]$, which is non-trivial in $H_2(W,\Z)$, is mapped trivially
to  $H_2(W, Y_1\cup Y_0,\Z)$. Finally, we have $H_2(W_2,Y_1\cup
Y_0,\Z)=\Z$ generated by $[\tilde \Sigma_0]$.
\end{Lem}

\proof   The results simply follow from the exact sequence in
homology. We have 
$$H_3(W_0, Y \cup Y_0,\Z)\cong H_1(W_0)=0, $$
hence we have
$$ 0 \to H_2(Y_0,\Z) \stackrel{\cong}{\longrightarrow} H_2(W_0,\Z)
\stackrel{0}{\longrightarrow} H_2(W_0, Y \cup Y_0,\Z)\longrightarrow
H_1(Y_0,\Z) \to 0. $$
Similarly, we have 
$$ 0 \to H_2(Y_0,\Z) \longrightarrow H_2(W,\Z)\longrightarrow H_2(W,
Y_1\cup Y_0,\Z) \longrightarrow
H_1(Y_0,\Z) \to 0 $$
which is of the form
$$ 0 \to \Z\to \Z\oplus \Z \stackrel{0\oplus I}{\longrightarrow}
\Z\oplus \Z \to \Z \to 0. $$

\endproof

Thus, we obtain the following result.

\begin{Lem}
We have an identification of the $\spinc$-structures on $W_0$ and
$W_1$ given by
$$ {\cal S}(W_0) =\{ \s_k \}_{k\in \Z}, $$
with
$$ \s_k \mapsto c_1(\det W_{\s_k}^+)= 2k \ e_0= 2k PD_{W_0}[\tilde
\Sigma_0], $$ 
with $[\tilde \Sigma_0]$ the generator of $H_2(W_0,Y\cup
Y_0,\Z)$. Similarly, we have  
$$ {\cal S}(W_1) =\{ \s_{\ell} \}_{\ell\in \Z} , $$
with
$$ \s_{\ell} \mapsto c_1(\det W_{\s_{\ell}}^+)=(2\ell +1) \ e_1=(2\ell 
+1) PD_{W_1}[\Sigma_1], $$  
with $[\Sigma_1]$ the generator of
$H_2(W_1,Y_1\cup Y,\Z)$.
Moreover, we have
$$ {\cal S}(W)=\{ \s_{\ell,k} \}_{\ell,k \in \Z}, $$
with
$$ \s_{\ell,k} \mapsto c_1(\det W_{\s_{\ell,k}}^+)= (2\ell+1) \ e + 2k \
e_0, $$
where we have $e=  PD_{W}(E)$, and $E$ is the image of the class $[\Sigma_1
-\Sigma_0]$ in $H_2(W,Y_1\cup Y,\Z)$. Similarly, we have
$$ {\cal S}(W_2)=\{ \s_k \}_{k\in \Z}. $$
\label{spin:cobord}
\end{Lem}

\proof
The result follows from the previous Lemma with the additional
observation that in $W_1$ we have $\Sigma_1$ with self intersection
$\Sigma_1^2 =-1$, hence the $\spinc$ structures have odd Chern class
$(2\ell +1)PD_{W_1}[\Sigma_1]$.
\endproof

In the case of $W_0$, $W_2$, and $W$, the $\spinc$
structure $\s_k$ or $\s_{\ell,k}$ restricts to the end $Y_0 \times
[T_0,\infty)$ of the cobordism to the pullback of the $\spinc$
structure $\s_k$ in ${\cal S}(W_0)$. All the $\spinc$ structures of
Lemma \ref{spin:cobord} restrict to the trivial $\spinc$ structure on
the cylindrical ends modeled on $Y$ or $Y_1$.

\subsection{Splitting of the cobordisms}

In the following, we shall introduce moduli spaces of Seiberg-Witten
equations on the cobordisms. Our purpose is to apply to the
moduli spaces on the cobordisms the same
techniques we developed in \cite{MW2}, in the study of the moduli spaces
of flowlines, that is, of monopoles on the trivial cobordisms. Thus,
it is convenient to consider the manifolds $W_1$ and $W_0$ endowed
with infinite cylindrical ends $Y_1 \times (-\infty,
-T_0]$ and $Y\times [T_0,\infty)$, and $Y_0 \times
[T_0,\infty)$ and $Y\times (-\infty, -T_0]$, respectively, with
metrics $g_Y +dt^2$ and $g_{Y_i}+dt^2$.  Moreover, we shall assume
that the 3-manifolds $Y_1$, $Y$, and $Y_0$ are endowed with metrics
with a long cylinder $T^2\times [-r,r]$, as specified in \cite{CMW}.

We can then think of the cobordisms as endowed with a metric which
restricts to the flat product metric on the region  
$T^2\times [-r,r]\times \R$. Moreover, we can identify in the
cobordisms $W_i$ a product region $V\times \R$, on the complement of a
tubular neighbourhood of the knot, where the cobordism is trivial. 
Thus, we can decompose the cobordisms $W_i$ as
\ba W_i= V\times \R \cup_{T^2\times \R} T^2\times [-r,r]\times \R
\cup_{T^2\times \R}  W_i(\nu(K)). \label{split:cobord} \na
The non-compact region $W_i(\nu(K))$ has the following property. There 
is a compact set ${\cal K}$ in $W_i$ such that the intersection
${\cal K} \cap W_i(\nu(K))$ is obtained by attaching a 2-handle 
$D\times D$ to the product $\nu(K)\times [-T_0, T_0]$, and, outside
of ${\cal K}$, the region ${\cal K}^c \cap W_i(\nu(K))$ consists of product
regions $\nu(K)\times [ T_0 ,\infty)$ and $\nu(K)\times (-\infty, -T_0]$,
and $T^2\times [r_0,r] \times [-T_0,T_0]$.

In the cobordism $W_i$ consider an interior point $x_i$
contained in the core disk of the 2-handle, $x_i \in D_i$. As in
\cite{BD}, we denote by $\hat W_i$ the 
punctured cobordism $\hat W_i = W_i \backslash \{ x_i \}$. Similarly,
we can consider the punctured manifold
$$ \hat W_i(\nu(K)) = W_i(\nu(K))\backslash \{ x_i \} . $$

In the manifolds $\hat W_i(\nu(K))$ we can identify a product region 
\ba {\cal V}=\nu(K)_{r_0}\times \R \cong D\times (D_i\backslash \{ x_i
\}). \label{prodW} \na 

This corresponds to endowing the manifold $\hat W_i$ with an extra
asymptotic end of the form $S^3 \times [0,\infty)$ at
the puncture. Thus, we identify the manifold $W_i$ with a
connected sum
$$ W_i = \hat W_i \# Q_i, $$
with a long cylindrical neck $S^3\times [-T(r),T(r)]$, and with $Q_i$ a
4-ball, as in \cite{BD}. 

\begin{figure}[ht]
\epsfig{file=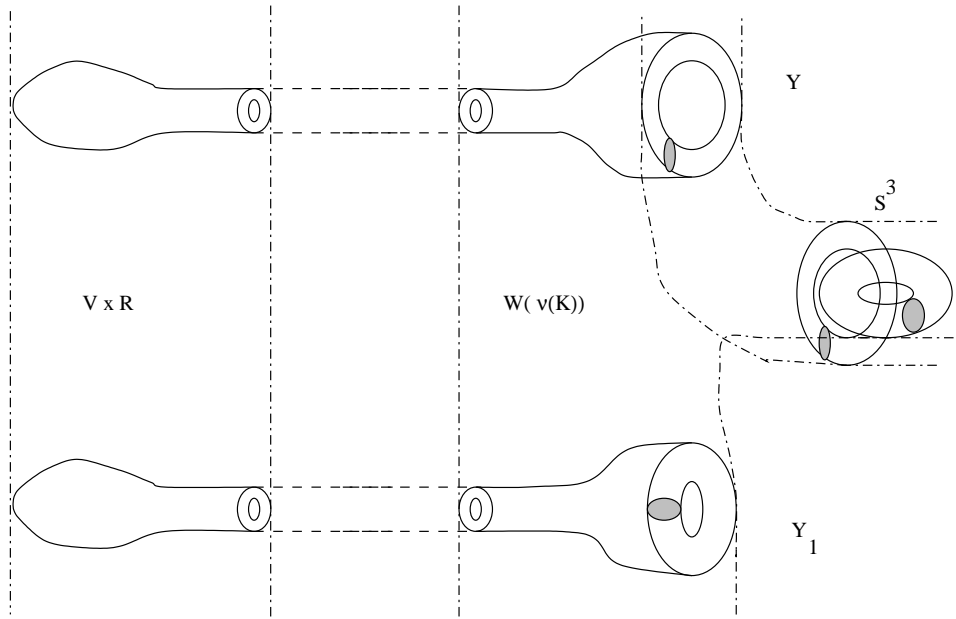,angle=270}
\caption{The decomposition of the punctured cobordism $\hat W_1$
\label{figIII2}} 
\end{figure}

Consider the sphere $S^3$ decomposed as the union of
two solid tori in the standard way,
$S^3=\nu \cup \tilde \nu$, with $\nu\cong \tilde\nu \cong D\times
S^1$. Then the product region ${\cal V}$ of (\ref{prodW}) 
in $W_i$ identifies the standard solid torus $\nu$ in $S^3$ with the
neighborhood 
$\nu(K)$ of the knot $K$ in $Y$, and, similarly, the other solid torus
$\tilde\nu$ in $S^3$ is identified with the 
tubular neighbourhood $\nu(K)$ in $Y_i$, after the surgery.
This is illustrated in Figure \ref{figIII2}.

\subsection{Metrics and perturbations on the cobordism}

The results of this subsection are based on the metric deformation Lemma
that Liviu Nicolaescu kindly communicated to us, \cite{Ni:pri}, and
that we enclosed in Part I \cite{CMW}.

\begin{Lem}
Let $A$ be an element in $SL(2,\Z)$. Suppose given $\epsilon >0$
sufficiently small. Consider the metric on $T^2$ given by
$$ g_0= A^* g, $$
where $g$ is the standard flat metric as before. There exists a
constant $\delta >0$ and a smooth path $g(s)$ of flat metrics on $T^2$ 
with the following properties:

(i) $g(s)\equiv \frac{1}{\delta^2} g_0$, for all $s \leq \epsilon$;

(ii) $g_1 = g(1)$ is a metric of the form
$$ g_1 = k_1 du^2 + k_2 dv^2 $$
for $k_i$ positive constants;

(iii) $g(s)=g_1$ for all $s\geq 1-\epsilon$;

(iv) The scalar curvature of the metric $\hat g= ds^2 + g(s)$ on $T^2
\times \R$ is non-negative.

(v) The metric $g_1$ can be extended to a metric inside the 
solid torus $\nu(K)$, which 
we still denote $g_1$, that has non-negative scalar curvature.

The constant $\delta$ is given by $\delta^2 = g_0(\partial_{u}, 
\partial_{u})$.
\label{Liviu:metrics}
\end{Lem}

Using this result, when we construct the manifold $Y_1$ from $Y$,
by removing a tubular neighborhood $\nu(K)$ and gluing it
back along $T^2$ with the matrix $A \in SL(2,\Z)$ prescribed by the
surgery, we can consider the same metric on the knot complement $V$,
with an end isometric to $T^2\times [0,\infty)$, with the metric $g+
ds^2$. On the other hand, on $\nu(K)$ we can consider the
metric $\delta^2(g(s) + ds^2)$, with $g(s)$ constructed as above, with the
parameterization chosen in such a way that we have 
$$ g(s)\equiv \frac{1}{\delta^2} A^*g $$
on the end $T^2\times [\delta^{-2}r_0,\infty)$ and $g(s)\equiv g_1$ 
near $T^2 \times \{ 0 \}$, extended to a positive scalar
curvature metric inside the solid torus.

With this choice of metrics, we still have the decomposition of the
moduli spaces of critical points of the $CSD$ functional, as proved in 
Part I, \cite{CMW}. Moreover,
this particular choice of metrics allows us to describe the choice of
metrics on the cobordisms.

\begin{figure}[ht]
\epsfig{file=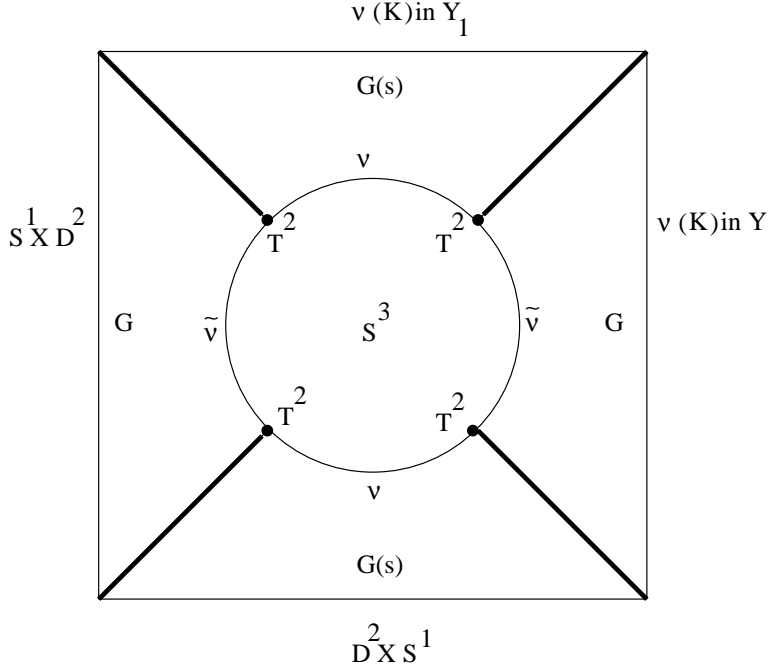,angle=270}
\caption{Product regions and metrics on the punctured 2-handle
\label{figIIIhandle}} 
\end{figure}

Consider first the trivial cobordism $Y_1
\times \R$. In the limit $r\to\infty$ this splits as $V\times \R$,
with an end of the form $T^2\times \R^+ \times \R$ with metric $g+
ds^2 + dt^2$, and $\nu(K)\times \R$, with an end of the form
$T^2\times \R^+ \times \R$ with the metric $\delta^2(g(s) + ds^2)+ dt^2$, as 
described above. 
Now consider the punctured cobordism $\hat W_1$. This contains a
product region $\nu\times \R$ which connects the solid torus $\nu(K)
\subset Y$ with a solid torus $\nu \subset S^3$ at the puncture.
On this product region we consider the metric $G:=g+ \delta^2 ds^2+ dt^2$, 
with $g$ the standard flat metric on $T^2$ extended to a non-negative scalar   
curvature metric inside the solid torus as described in \cite{CMW}. 
The other product region $\tilde \nu \times \R$ connecting the solid
torus $\nu(K)\subset Y_1$ with a solid torus $\tilde \nu$ inside $S^3$,
is glued in the punctured 2-handle, along a region $T^2\times \R$,
with framing one. Thus, on this product region we can consider the metric
$G(s):=\delta^2(g(s)+ds^2)+dt^2$ described above.
These regions are illustrated schematically in
Figure \ref{figIIIhandle}, with a lower dimensional picture of the
punctured handle.

Our purpose is to define chain maps between the Floer complexes
of the 3-manifolds using the Seiberg--Witten equations on the
cobordisms, and to adapt the techniques of \cite{MW2} to analyze these 
chain maps, by understanding their asymptotic limits under the
splitting of the cobordisms illustrated in the previous section. Thus, 
we need to introduce a suitable perturbation of the Seiberg--Witten
equations on $W_i$ which is compatible with the perturbations of the
Chern--Simons--Dirac functional on the manifolds $Y_i$ and $Y$,
described in \cite{CMW}.

Recall
that on $Y_1$ we have perturbed flow equations of the form
\ba
\left\{ \begin{array}{lll}
\displaystyle{\frac{\d A}{\d t} }&=& -*F_A + \sigma (\psi, \psi) + 
\sum_{j=1}^{N} \displaystyle{\frac{\d U_1}{\d \tau_j}} \mu_j^{(1)} \\[2mm]
\displaystyle{\frac{\d \psi}{\d t}}& =& -\dirac_A \psi 
- \sum_{i=1}^{K} \displaystyle{\frac{\d V_1}{\d \zeta_i}}\nu_i^{(1)}.\psi,
\end{array}\right.
\label{CSD-F}
\na
where $(U_1,V_1)$ is a pair of functions in
the class ${\cal P}_\delta$ described in
\cite{CMW}, that is, it becomes exponentially small along the cylinder 
$T^2\times [-r,r]$ inside $Y_1$, and is exponentially small on the
solid torus $\nu(K)\subset Y_1$. All the notation we use here follows
\cite{CMW}.   

The equations (\ref{CSD-F}) can be written equivalently as the four
dimensional equations on $Y_1\times \R$,
\ba \left\{\begin{array}{l}
F_{\A}^+ = \tau(\Psi,\Psi) + P_1^{Y_1}(\A,\Psi) \\[2mm]
D_{\A}\Psi=P_2^{Y_1}(\A,\Psi), \end{array}\right. \label{4SWP}
\na
with
\ba
P_1^{Y_1}(\A,\Psi)=\sum_{j=1}^{N} \displaystyle{\frac{\d U_1}{\d \tau_j}}
(\mu_j^{(1)}\wedge dt + *_3 \mu_j^{(1)}), \label{P1} \na
\ba
P_2^{Y_1}(\A,\Psi)=\sum_{i=1}^{K} \displaystyle{\frac{\d V_1}{\d
\zeta_i}}\nu_i^{(1)}.\Psi. \label{P2} \na 

As in \cite{CMW}, we consider similarly perturbed equations on $Y$,
with the additional surgery perturbation
on the solid torus $\nu(K)\subset Y$. Thus, on $Y$ we have 
flow equations
\ba
\left\{ \begin{array}{lll}
\displaystyle{\frac{\d A}{\d t} }&=& -*F_A + \sigma (\psi, \psi) + 
\sum_{j=1}^{N} \displaystyle{\frac{\d U}{\d \tau_j}} \mu_j +
f'(T_{A(t)}) *\mu \\[2mm] 
\displaystyle{\frac{\d \psi}{\d t}}& =& -\dirac_A \psi 
- \sum_{i=1}^{K} \displaystyle{\frac{\d V}{\d \zeta_i}}\nu_i.\psi,
\end{array}\right.
\label{CSD-F-Y}
\na
with 
\[
T_A(z)= -i\int_{\{z\in D^2\} \times S^1} (A-A_0).
\]
Here the function $f'$, which depends on the choice of a small parameter
$0<\epsilon \leq \epsilon_0$, is constructed as in \cite{CMW}, $\mu$ is
a compactly supported 2-form, and $(U,V)$ is a perturbation in the
class ${\cal P}_\delta$ on $Y$.

The equation (\ref{CSD-F-Y}) can be written equivalently as 
\ba \left\{ \begin{array}{l}
F_{\A}^+ = \tau(\Psi,\Psi) + P_1^Y(\A,\Psi) + f'(T_{\A})\mu^+\\[2mm]
D_{\A}\Psi=P_2^Y(\A,\Psi), \end{array}\right. \label{4SWP:surg}
\na
where $\mu^+$ is the self-dual part of the pullback of $\mu$ along
the projection $\nu(K)\times \R \to \nu(K)$.

We discuss a choice of perturbation on $W_i$ which behaves nicely
under the splitting of the cobordisms and restricts to the prescribed
perturbations on the asymptotic ends. Consider the 
case of $W_1$. The case of $W_0$ is analogous. 

The manifold $W_1$ has a cylindrical end $Y\times [T_0,\infty)$ and a
cylindrical end $Y_1 \times (-\infty, -T_0]$. 
Consider a cutoff function $\chi(t)$ which is  
$\chi(t)\equiv 1$ for $t \geq T_0$ and $\chi(t)\equiv 0$ for  $t
\leq T_0 -1$. In the following we shall use the notation
$\hat\chi(t)=\chi(-t)$. 
Consider, on the cylindrical ends $Y\times [T_0,\infty)$  and $Y_1
\times (-\infty, -T_0]$ of $W_1$, the
equations (\ref{4SWP}) and (\ref{4SWP:surg}), respectively.
Now consider the manifold $V\times \R$, endowed with an infinite end of
the form $T^2\times [0,\infty)\times \R$. Notice that on this region
inside the manifold $V\times \R$
it makes sense to define the temporal gauge condition for pairs
$(\A,\Psi)$. We denote by $(A(t),\psi(t))$ a temporal gauge
representative of the gauge class of $(\A,\Psi)$.

On $V$ we can also
introduce perturbations in the class ${\cal P}_\delta$, as described in
Section 3 of\cite{CMW}. 

This gives a choice of perturbation on  $V_r \times \R$ inside
$W_1$, for large $r\geq r_0$. Namely, we consider on $V_r \times \R$
the equations
\ba \left\{ \begin{array}{ll}
F_{\A}^+ = & \tau(\Psi,\Psi) + (1-\hat\chi -\chi) P_1^V(\A,\Psi)
\\[2mm] & + \hat\chi P_1^{Y_1}(\A,\Psi) + \chi P_1^Y(\A,\Psi) \\[2mm] 
D_{\A}\Psi= & (1-\hat\chi-\chi) P_2^V(\A,\Psi) \\[2mm]
& + \hat\chi P_2^{Y_1}(\A,\Psi) + \chi
P_2^Y(\A,\Psi). \end{array}\right. 
\label{4SWP-V} 
\na 
This matches the perturbations (\ref{4SWP}) and (\ref{4SWP:surg}) on
$V\times [T_0,\infty)$ and $V\times (-\infty, -T_0]$ in $W_1$.

On the product region ${\cal V}$ of (\ref{prodW}),
which connects $\nu(K) \subset Y$ to the solid torus $\nu\subset S^3$,
we consider the Seiberg--Witten equations with the surgery
perturbation, namely the equations
\ba \left\{
\begin{array}{l}
F_{\A}^+ =  \tau(\Psi,\Psi) + f'(T_{\A})\mu^+  
+ \chi P_1^Y(\A,\Psi)\\[2mm] 
D_{\A}\Psi =  \chi
P_2^Y(\A,\Psi). \end{array}\right. \label{4:eq:surg} \na 

On the  product region $\tilde {\cal V}$ inside $\hat
W_1(\nu(K))$ that connects the solid torus $\tilde \nu$ in $S^3$ to
$\nu(K)$ in $Y_1$, we consider equations
\ba \left\{
\begin{array}{l} 
F_{\A}^+ = \tau(\Psi,\Psi)+ \hat\chi P_1^{Y_1}(\A,\Psi)  \\[2mm]
D_{\A}\Psi = \hat\chi P_2^{Y_1}(\A,\Psi). \end{array}\right. \label{4:eq}
\na

In the case of $W_0$, if we consider $Y_0$ with the trivial
$\spinc$-structure, the perturbation along the cylindrical end
$Y_0\times [T_0,\infty)$ includes a perturbation $\rho_0$ with
$[*\rho_0]=\eta PD_{Y_0}(m)$, with respect to the $*$-operator
on $Y_0$, with $\eta>0$ as in \cite{CMW} and $m$ the generator of the
first homology group of $Y_0$. 

Throughout the paper, when we consider finite energy solutions of the
Seiberg--Witten equations on the cobordisms $W_i$, we will mean finite
energy solutions of the Seiberg--Witten equations on the punctured
cobordisms $\hat W_i$, as in Definition \ref{fin:en:def} below,
with the perturbations introduced here, and with a
removable singularity at the puncture $x_i$, that is, such that they
extend to solutions on $W_i$.

In the following, we shall analyze the behavior of finite energy
solutions on $W_i$, when stretching $r\to \infty$ in the cobordism
\ba \hat W_i(r) = V_r\times \R \cup {\cal V}(r) \cup \tilde {\cal
V}(r), \label{W:parts} \na 
with ${\cal V}(r)\cong \nu_r \times \R$ and $ \tilde {\cal V}(r)\cong
\tilde \nu_r\times \R$, the product regions
inside  $\hat W_i(\nu(K))_r$.

In particular, with the choice of perturbation discussed
here, a finite energy solution on $W_1$ will have an asymptotic value
$$ a_1=[A_1,\psi_1] \in {\cal M}_{Y_1} $$
on the $Y_1$ end,
and
$$ a=[A,\psi] \in {\cal M}_{Y,\mu}, $$
on the other end of the cobordism, with the moduli spaces of solutions
of the perturbed equations on $Y$ and $Y_1$ as in \cite{CMW}.
Similarly, a finite energy solution on the cobordism $W_0$ will have
asymptotic values
$$ a=[A,\psi] \in {\cal M}_{Y,\mu}, $$
$$ a_0=[A_0,\psi_0] \in {\cal M}_{Y_0}(\s_k). $$
Thus we can define moduli spaces $\M^{W_1}_{\s_\ell}(a_1,a)$ and
$\M^{W_0}_{\s_k}(a,a_0)$ as follows.

\subsection{Finite energy monopoles and virtual dimension}

We consider finite energy solutions of the perturbed equations
introduced above, on the punctured cobordisms $\hat
W_i$, with a removable singularity at the puncture $x_i$, that is,
such that they extend to solutions on $W_i$. More
precisely, we have the following.

\begin{Def} \label{fin:en:def}
Consider the manifold $\hat W_i$, as a complete Riemannian manifold
with infinite cylindrical ends. Consider a fixed choice of the
$\spinc$ structure. Write $C$ for any of the cylindrical ends of the
manifold $\hat W_i$.
A solution $(\A,\Psi)$ of the Seiberg--Witten equations
on the manifold $\hat W_i$ is of finite energy, with a removable
singularity at the point $x_i$, iff on any of the
cylindrical ends $C$ the solutions $(\A,\Psi)$ in a temporal gauge 
satisfies the condition 
$$ \int_C | \partial_t A(t) |^2 + |\partial_t \psi(t)|^2  <\infty. $$ 
\end{Def}

The analysis in Section 3 of \cite{MW} of the asymptotics of finite energy
solutions on a trivial cobordism $Y\times \R$ carries over to the
present case and shows that finite energy solutions on the manifolds
$\hat W_i$ decay along the cylindrical ends to asymptotic values $(A,\psi)$
satisfying the 3-dimensional 
Seiberg--Witten equations on the boundary 3-manifolds.
Moreover, if the asymptotic value is an irreducible critical point,
then the rate of decay is exponential, with the exponent determined
by the first eigenvalue of the Hessian at the critical point.

Thus, we define configuration spaces ${\cal A}_{k,\delta}(W_i)$ and
the group of gauge transformations ${\cal G}_{k+1,\delta}(W_i)$ and
consider finite energy solutions in the quotient space. 
The configuration space ${\cal A}_{k,\delta}(W_i)$ consists of pairs
$(\A,\Psi)$ on $\hat W_i$ that are of finite energy, with a removable
singularity at $x_i$, and with a rate of decay with exponent $\delta$
along the cylindrical ends modeled on $Y$ and $Y_1$. Since the
asymptotic limit along the end $S^3\times [0,\infty)$ at the puncture
$x_i$ is the unique reducible solution $\theta_{S^3}$ on $S^3$, we can
define moduli spaces
$$ \M^{W_1}_{\s_\ell}(a_1,a) \ \ \hbox{ and } \ \ \M^{W_0}_{\s_k}(a,a_0) $$
of solutions in ${\cal A}_{k,\delta}(W_i)$ modulo gauge action, which
depend only on the asymptotic limits at the two ends of $W_i$ and on
the $\spinc$ structure on $W_i$.
We can estimate the virtual dimension of these moduli spaces of solutions.

\begin{Lem}
Consider the linearization ${\cal D}_{(\A,\Psi)}$ of the 
Seiberg--Witten equations at a solution $(\A,\Psi)$
on $(W_1, \s_\ell)$, with asymptotic values $a_1=[A_1,\psi_1]$ and $a=[A,
\psi]$ on the boundary three-manifolds $Y_1$ and $Y$. 
The index of ${\cal D}_{(\A,\Psi)}$ is given by
$$ Ind({\cal D}_{(\A,\Psi)})=\frac{1}{4}(c_1(\det(\s_\ell))^2 -2\chi(W_1)
-3\sigma(W_1)) + \frac{\rho}{2}(Y_1,a_1) - \frac{\rho}{2}(Y,a), $$
where the last two summands are the APS $\rho$-invariants \cite{APS}
of the extended Hessian operators $H_{(A_1,\psi_1)}$ and
$H_{(A,\psi)}$. The grading analyzed in \cite{CMW} satisfies
$$ \deg_{Y_1} (a_1) =  \frac{1}{8 \pi^2} CSD_{Y_1}(a_1) +
\frac{\rho}{2}(Y_1,a_1), $$ 
up to redefining the $CSD$ functional by a global additive constant.
\end{Lem}

\proof
In the case of a solution $(\A,\Psi)$ on the trivial cobordism
$Y\times I$,  
with asymptotic values $a=[A_a,\psi_a]$ and $b=[A_b,
\psi_b]$ on $Y$, we have
$$ \deg_Y (a) -\deg_Y (b) = Ind({\cal D}_{(\A,\Psi)})= 
\frac{c_1(\det(\s))^2}{4} +\frac{\rho}{2}(Y,a) - \frac{\rho}{2}(Y,b)$$
$$ = \frac{-1}{16 \pi^2}\int_{Y\times\R} F_{\A}\wedge F_{\A}
+\frac{\rho}{2}(Y,a) - \frac{\rho}{2}(Y,b) $$
$$ = \frac{-1}{8 \pi^2}(CSD(A_b,\psi_b)- CSD(A_a,\psi_a)) +
\frac{\rho}{2}(Y,a) -
\frac{\rho}{2}(Y,b). $$ 

\endproof

See the energy estimates in \cite{MW2} (cf. \cite{APS}, \cite{MMR}).
This virtual dimension of Seiberg--Witten moduli spaces on
four-manifolds with boundary is computed explicitly in
\cite{Nic}, in the case of Seifert fibered spaces, where a particular
choice of metric makes it possible to compute $\rho(Y,a)$ explicitly.

By the results of \cite{CMW}, the choice of the grading on $Y_1$
determines uniquely the grading $\deg_{Y,\mu}$ on $Y$, up to changing
the functional $CSD$ by a global additive constant. The grading
$\deg_{Y,\mu}$ also  determines \cite{CMW} the choice of 
the grading $\deg_{Y_0}$ on $Y_0$. For the properties of the Floer
complex on $Y_0$ see also \cite{MW4}.

We define the expression
$$ \iota(\s_\ell,W_1,a_1,a)=\frac{1}{4}(c_1(\det(\s_\ell))^2 -2\chi(W_1)
-3\sigma(W_1)) + \frac{\rho}{2}(Y_1,a_1) - \frac{\rho}{2}(Y,a) $$
in the case of $W_1$. Similarly, in the case of  $W_0$,
we define
$$ \iota(\s_k,W_0,a,a_0)=\frac{1}{4}(\int_{W_0}c_1(\det(\s_k))^2
-L(\hat\nabla^0))-\frac{1}{2}(\chi(W_0) 
+\sigma(W_0)) $$
$$ + \frac{\rho}{2}(Y,a) - \frac{\rho}{2}(Y_0,a_0), $$
where $\hat\nabla^0$ is the metric compatible connection, which along
the end $Y_0 \times [T_0,\infty)$ has the form $dt\otimes \partial_t + 
\nabla^0$.(This is analogous to the case of \cite{Nic}, but in our
setting the connection $\nabla^0$ is simply the metric connection with 
no adiabatic rescaling.) The expressions
$$ \iota(\s_\ell,W_1,a_1,a) \ \ \hbox{and} \ \ \iota(\s_k,W_0,a,a_0)$$
compute the virtual dimensions of the moduli spaces 
$$ \M^{W_1}_{\s_\ell}(a_1,a) \ \ \hbox{ and } \ \ \M^{W_0}_{\s_k}(a,a_0). $$
Notice that, in the case of $W_0$,  the
$\spinc$-structure may be non-trivial along the end $Y_0 \times
[T_0,\infty)$. This has the important consequence that the virtual dimension
$\iota(\s_k,W_0,a,a_0)$
is only defined modulo the integer $d(\s_k)$, with $d(\s_k)$ satisfying
$$ c_1(\det(\s_k))(H_2(Y_0,\Z))=d(\s_k)\Z, $$
that is, in this case, $d(\s_k)=2k$. This ambiguity corresponds to
different components of the moduli space of solutions of the
Seiberg-Witten equations on $W_0$, with different energies. The
minimal energy component corresponds to the minimal non-negative value 
of $\iota(\s_k,W_0,a,a_0)$.

\subsection{Compactification and invariants}

Under a generic choice of the perturbation, we can
assume that all the moduli spaces $\M^{W_1}_{\s_{\ell}}(a_1,a)$ and
$\M^{W_0}_{\s_k}(a,a_0)$ are cut out transversely, of dimension
$\iota(\s_{\ell},W_1,a_1,a)$ and $\iota(\s_k,W_0,a,a_0)$,
respectively. 
Unless otherwise stated, when we write $\M^{W_0}_{\s_k}(a,a_0)$, we
only consider the component of minimal energy, with
dimension $\iota(\s_k,W_0,a,a_0)$. 

The following description of the compactification of
$\M^{W_1}_{\s_{\ell}}(a_1,a)$ and $\M^{W_0}_{\s_k}(a,a_0)$ follows
from the main gluing theorem proved in \cite{MW}, together with the 
results of the previous subsections.

\begin{Pro}
Suppose given a non-empty moduli space $\M^{W_0}_{\s_{k}}(a,a_0)$
of dimension
$$ \iota(\s_{k},W_0,a,a_0)= n >0, $$
with $a\in \M_{Y,\mu}^*$.
Then $\M^{W_0}_{\s_{k}}(a,a_0)$ admits a compactification to a manifold
with corners, where the  
codimension $1$ boundary strata consist of
\ba \label{strata0}
\begin{array}{c}
\bigcup_{c\in \M_{Y,\mu}^*}  \hat \M_{Y,\mu}(a,c)\times  \M^{W_0}_{\s_{k}}(c,
a_0) \\[2mm]
\bigcup_{c_0\in \M_{Y_0}(\s_k)} \M^{W_0}_{\s_{k}}(a,c_0) \times 
\hat \M_{Y_0,\s_k}^{(0)}(c_0,a_0),\end{array}
\na 
and $\hat \M_{Y_0,\s_k}^{(0)}(c_0,a_0)$ is the minimal energy component
of the moduli space $\hat \M_{Y_0,\s_k}(c_0,a_0)$, as discussed in
\cite{MW4}. If the reducible point $c=\theta$ satisfies
$$ \deg_{Y,\mu}(a)>\deg_{Y,\mu}(\theta), $$
then we also have an extra component  
in (\ref{strata0}) of the form
\ba \hat \M_{Y,\mu}(a,\theta)\times U(1) \times
\M^{W_0}_{\s_{k}}(\theta, a_0)  \label{redglue0} \na
with a $U(1)$ gluing parameter.

We have a similar compactification of $\M^{W_1}_{\s_\ell}(a_1,a)$
of dimension
$$ \iota(\ell,W_1,a_1,a)=n >0, $$
for $a_1\in\M_{Y_1}^*$ and $a\in\M_{Y,\mu}^*$,
with codimension $1$ boundary strata of the form
\ba \label{strata1}
\begin{array}{c}
\bigcup_{c\in \M_{Y,\mu}^*}\M^{W_1}_{\s_\ell}(a_1,c)\times \hat \M_{Y,\mu}(c,a)
\\[2mm]
\bigcup_{c_1 \in \M_{Y_1}^*} \hat\M_{Y_1}(a_1,c_1)\times
\M^{W_1}_{\s_\ell}(c_1,a), \end{array} \na
and with extra components
\ba \begin{array}{c} \hat \M_{Y_1}(a_1,\theta_1)\times U(1) \times
\M^{W_1}_{\s_\ell}(\theta_1,a) \\[2mm]
\M^{W_1}_{\s_\ell}(a_1,\theta)\times U(1) \times\hat
\M_{Y,\mu}(\theta,a), 
\end{array} \label{redglue1} \na
when splitting through the reducibles. 
\label{compactif}
\end{Pro}

The Proposition follows from the main gluing theorem of \cite{MW}. 
The fact that only the minimal energy component, 
among the components of $\hat
\M_{Y_0}(\lambda c_0, a_0)$, occurs in the
compactification of $\M^{W_0}_{\s_{k}}(c,a_0)$ is explained in the
following, in the proof  of Lemma \ref{chain:homom}.

For later use, we also need the following.

\begin{Cor}
If we have $ \iota(\s_\ell,W_1,a_1,a)=1$, or $
\iota(\s_{k},W_0,a,a_0)=1$, with $a_1\in\M_{Y_1}^*$,
$a\in\M_{Y,\mu}^*$, and $a_0\in \M_{Y_0}(\s_k)$, then the boundary
strata are given by 
(\ref{strata1}) and (\ref{strata0}), as in the compactification of
Proposition \ref{compactif}, with $c_1$ and $c$ irreducible.
\label{irredonly}
\end{Cor} 

\proof We need to show that the components (\ref{redglue0}) and
(\ref{redglue1}) do not occur in the compactification. This follows
by dimensional arguments. 
\endproof

We can define numerical invariants associated to the zero-dimensional
components of the moduli spaces 
$\M^{W_1}_{\s_{\ell}}(a_1,a)$, and $\M^{W_0}_{\s_k}(a,a_0)$, with
$a_1\in \M_{Y_1}^*$,  
$a\in \M_{Y,\mu}^*$, and $a_0\in \cup_k \M_{Y_0}(\s_k)$. 

Recall that these moduli spaces come endowed with an orientation,
given by the trivialization of the determinant line bundle of the
linearization of the Seiberg--Witten equations. Some care is needed in 
defining the orientation in the case of non-compact 4-manifolds with
infinite cylindrical ends. The necessary details can be found in
\cite{Nic}. The orientation is compatible with the compactification of 
Proposition \ref{compactif}. Throughout this discussion we shall always
assume that the 
perturbations are chosen so that all the moduli spaces are cut out
transversely by the equations.

According to Proposition \ref{compactif} of the 
previous section, in the case of $\iota(\s_{\ell},W_1,a_1,a)=0$, or
$\iota(\s_k,W_0,a,a_0)=0$, the moduli spaces
$\M^{W_1}_{\s_{\ell}}(a_1,a)$ and $\M^{W_0}_{\s_k}(a,a_0)$ consist of
a finite set of points with an attached sign given by the orientation.
Thus, we can define
$$ N^{W_1}_{\s_{\ell}}(a_1,a) \ \ \hbox{ or } \ \ N^{W_0}_{\s_k}(a,a_0) $$
as the algebraic sum of the points in $\M^{W_1}_{\s_{\ell}}(a_1,a)$
and $\M^{W_0}_{\s_k}(a,a_0)$, respectively. If we have 
$\iota(\s_{\ell},W_1,a_1,a)<0$ and either $a$ or $a_1$ is
irreducible, or if we have $\iota(\s_k,W_0,a,a_0)<0$, then
the corresponding moduli space is generically empty, so we just set the 
corresponding invariant equal to zero.

\section{Geometric limits}

In this section we describe the geometric limits of finite energy
solutions on $W_i$ when stretching $r\to\infty$ inside the cobordism
$$ \hat W_i(r)= V_r\times \R \cup {\cal V}(r) \cup \tilde {\cal V}(r). $$
The analysis is very similar to the analysis of the geometric limits
of flow lines in \cite{MW2}. 

We give the general description of the geometric limits. 
We describe the case of the cobordism $W_1$. Simple modifications 
adapt the argument to $W_0$.
We shall omit here the parts of the
argument which are completely analogous to the case discussed in
\cite{MW2}. 
We assume that the metric is chosen as discussed previously.

\begin{Pro}
Consider a family $[\A_r,\Psi_r]$ of finite energy
solutions of the Seiberg--Witten equations on the cobordism $W_1(r)$,
with the choices of perturbation as discussed previously. Assume that
the $[\A_r,\Psi_r]$ have asymptotic values $a_1=[A_1,\psi_1]$ and $a=[A,
\psi]$ on the ends modeled on  $Y_1$ and $Y$. By \cite{CMW}, these
asymptotic limits can be written as 
$$ (A_1,\psi_1)=(A_1',\psi_1')\#_{a^-} (a^-,0), $$
and
$$ (A,\psi)=(A',\psi')\#_{a^+} (a^+,0), $$
for large $r\geq r_0$.

Let $\vartheta_1$ in $\chi_0(T^2,Y_1)$ denote the intersection point
between the lines $\{ v=0 \}$ and $\{ v-u=1 \}$. 

We have the following types of geometric limits of $[\A_r,\Psi_r]$, as
$r\to\infty$. 

(a) A finite energy solution $(\A',\Psi')$ of the perturbed equations
(\ref{4SWP}) on $V\times \R$. In radial gauge, this solution decays in
the radial direction to a flat connection $a_\infty '$ on $T^2$, 
$$ a_\infty ' \in \chi_0(T^2,V). $$ 
In a temporal gauge on $V\times \R$, $(\A',\Psi')$ converges to
elements $[A,\psi]$ and 
$[\tilde A,\tilde \psi]$ in $\M_V^*$ as $t\to\pm\infty$, with
$$ \partial_\infty [A,\psi] =a_\infty ' =\partial_\infty [\tilde
A,\tilde \psi]. $$

(b) Non-uniform limits on $V\times \R$ given by paths $[A(t),\psi(t)]$
in $\M_V^*$ 
connecting $[A,\psi]$ to $[A_1',\psi_1']$ and $[\tilde A,\tilde \psi]$
to $[A',\psi']$, and by a function $a_V : D^+ \to \chi_0(T^2,V)$,
holomorphic on some 
neighbourhood of the half disk $D^+$, which agrees on the subset
 $\theta\in \{ -\pi/2, \pi/2 \}$ of the boundary of $D^+$  with  the
asymptotic values $\partial_\infty[A(t),\psi(t)]=: a(t)$ of the paths in
$\M^*_V$. 

(c) A flat connection $a_\infty ''$ on $T^2$, obtained as a finite
energy solution of the equation
$$ F_{\A}^+ = f'(T_{\A}) \mu^+, $$
on the product region ${\cal V}=\nu(K)\times \R$ with an
infinite end $T^2\times [0,\infty)\times \R$.

(d) A flat connection  $\tilde a_\infty ''$ on $T^2$, obtained as a
finite energy solution of the unperturbed equations (\ref{4:eq}) on
the region
$$ \hat W_1(\nu(K))\backslash {\cal V} $$
with an infinite end $T^2\times [0,\infty)\times \R$.

(e) Non-uniform limits on the ends of $W_1(\nu(K))$ given by a path
$\tilde a''(\tau)$ in  
$$ \M^{red}_{\nu(K)}=\{ u-v=1 \},   $$ 
for $\tau\in [-1,1]$, with $\tilde a''(-1)=a^-$, $\tilde a''(0)=\tilde
a_\infty''$, 
and $\tilde a''(1)=\vartheta_1$, and a 
path $a''(\tau)$ in the perturbed
$$ \M^{red}_{\nu(K),\mu}=\{ v=f'(u) \}, $$
for $\tau\in [-1, 1]$, with $a''(-1)=\vartheta_1$, $a''(0)=a_\infty''$
and $a''(1)=a^+$.

(f) We also have a map
$a_{\nu(K)}: D^+_\epsilon \to H^1(T^2,\R)$,
holomorphic on some 
neighbourhood of the domain 
$$  D^+_\epsilon =\{ \rho \in [\log\epsilon,0], \theta \in
[-\pi/2,\pi/2] \}. $$ 
Upon identifying $\tau=e^\rho \sin\theta$, this map agrees  with  the
path 
$\tilde a''(\tau)$, 
for $\tau\in [-1, -\epsilon]$ and with $a''(\tau)$ for $\tau \in
[\epsilon,1]$, on the  
subset
$$ \theta\in \{ -\pi/2, \pi/2 \} \ \ \hbox{and} \ \ \rho \in
[\log\epsilon,0] $$ 
of the boundary of $D^+_\epsilon$. 

(g) A ``thin holomorphic triangle'' $\Delta(a_\infty'', \tilde
a_\infty'', \vartheta_1)$ in the character variety
$\chi_0(T^2,V)$ (cf. \cite{BD} pg.234), with vertices
$\{ a_\infty'', \tilde a_\infty'', \vartheta_1 \}$ and with two sides
along the lines $\{ v-u=1 \}$ and $\{ v= f'(u) \}$, with the
parameterization $\tilde a''(\tau)$ for $\tau \in [0,1]$ and
$a''(\tau)$ for $\tau \in [-1,0]$. 

(h) A limit on compact sets in the region $T^2\times [-r,r]\times \R$,
given by a flat connection $a_\infty$ on $T^2$, and a 
non-uniform limit after rescaling,
given by a map
$\hat a: D \to H^1(T^2,i\R)$, 
holomorphic up to the boundary, matching the
values of $a_V$ and $a_{\nu(K)}$, as described in \cite{MW2}.

\label{geom:lim}
\end{Pro}

\proof
The result follows from the analysis of the convergence of flow  lines
in \cite{MW2}. The limits (a), (b), (f) and (h) are derived exactly
as  the analogous cases  in \cite{MW2}.
The case (c) describes the limit of the solutions $(\A_r,\Psi_r)$,
uniformly on compact sets, in the product region ${\cal V}(r)$ as
$r\to\infty$. Up to gauge transformations, and up to passing to a
subsequence, the solutions $(\A_r,\Psi_r)$ converge smoothly on
compact sets in ${\cal V}(r)$ to a finite energy solution of the perturbed
abelian ASD equation 
$$ F_{\A}^+ = f'(T_{\A})\mu^+ $$
on the strip ${\cal V}$ with an infinite end
$T^2\times [0,\infty)\times \R$. By the analysis of \cite{MW2}, this
is (up to gauge) a constant flat connection $a_\infty ''$ on $T^2$,
with holonomies satisfying $v=f'(u)$.
Similarly, the case (d) describes the uniform convergence in $\tilde
{\cal V}(r)$. The solutions $(\A_r,\Psi_r)$ 
converge smoothly on compact sets in $\tilde {\cal V}(r)$
to a finite energy solution of the unperturbed abelian ASD equation on
the region $\tilde {\cal V}$ with an infinite end
$T^2\times [0,\infty)\times \R$. Again, by the analysis of \cite{MW2}
we see that this is up to gauge a constant flat connection $\tilde
a_\infty ''$ on $T^2$ contained in the line $\{ v-u=1 \}$ in
$\chi_0(T^2,Y_1)$. Case (e) describes the non-uniform limits in the
regions $\tilde {\cal V}(r)$ and ${\cal V}(r)$, away from compact
sets, after suitable rescaling as described in \cite{MW2}. We have
$\tilde {\cal V} \cong \tilde \nu \times \R$ and ${\cal V} \cong \nu
\times \R$, connecting the two solid
tori $\nu$ and $\tilde \nu$ in the standard Heegaard splitting of
$S^3$ to the solid tori $\nu(K)$ in $Y_1$ and $Y$, respectively. Thus,
we can adapt the analysis used in \cite{MW2} for
non-uniform limits on $\nu(K)\times \R$. After a suitable 
rescaling, we resulting non-uniform limits in the region $\tilde {\cal
V}$ consist of a path $\tilde a''(\tau)$ along the line $\{ u-v=1 \}$ in
$\chi_0(T^2,V)$ and a map $\tilde a_\nu : D^+ \to H^1(T^2,i\R)$,
holomorphic in a 
neighbourhood of $D^+$, which agrees with $\tilde a''$ along the
subset 
$$ \theta\in \{ -\pi/2, \pi/2 \} \ \ \hbox{and} \ \ \rho \in
(-\infty,0] $$ 
of the boundary of $D^+$. Similarly, the non-uniform limits on ${\cal
V}$, after rescaling, consist of a path $a''(\tau)$ along the curve $v=f'(u)$ in
$\chi_0(T^2,V)$ and a map $a_\nu : D^+ \to H^1(T^2,i\R)$, 
holomorphic in a
neighbourhood of $D^+$, which agrees with $a''(\tau)$  along the
subset 
$$ \theta\in \{ -\pi/2, \pi/2 \} \ \ \hbox{and} \ \ \rho \in
(-\infty,0] $$ 
of the boundary of $D^+$.

The thin holomorphic triangle of Case (g) is obtained by the overlap
of these two regions, with the vertex $\vartheta_1$
at the flat connection on $T^2=\partial \nu=\partial \tilde \nu$ which
extends to both sides of the standard Heegaard splitting of $S^3$ to
give the unique reducible solution on $S^3$. That is, we have
$\theta_{S^3}|_{T^2}= \vartheta_1$.
Notice that, in general, the flat connections $a_\infty '$, $a_\infty
''$, $\tilde a_\infty''$, $a_\infty$, $a^+$, and $a^-$ on $T^2$ are
all distinct.

\endproof

\subsection{The holomorphic triangles}

This subsection contains some observations on the holomorphic
triangles that appear among the geometric limits of solutions on the
cobordisms, as discussed in Proposition \ref{geom:lim}. A better
understanding of these triangles will be very useful in analyzing the
different properties of the coefficients of the chain maps defined by
the cobordisms $W_1$, $W_0$, and $W$.

In the following, let $\epsilon$ be the parameter used in the
construction of the surgery perturbation. 

\begin{Def}
Consider the unique holomorphic triangle $\Delta^\epsilon$ with
vertices $\{ a^-, \vartheta_1, a^+ \}$ and sides along the curves
$\ell_1= \{ v-u=1 \}$, $\ell_\mu=\{ v=f'(u) \}$ and
$\ell=\partial_\infty(\M_V^*)$, defined by the geometric limits 
$a_V$, $a_{\nu(K)}$, $\hat a$, and $\Delta(\tilde
a_\infty'',\vartheta_1,a_\infty'')$ of Proposition \ref{geom:lim}.
We say that the triangle $\Delta^\epsilon$ is {\em
degenerate} if the holomorphic map obtained as a limit of the
triangles $\Delta^\epsilon$, as $\epsilon \to 0$ is a 
disk $\Delta$ with boundary along arcs of $\ell_1$ and
$\ell$ connecting the vertices $a^-$ and $a^+$. We say that
$\Delta^\epsilon$ is {\em non-degenerate} if the holomorphic map
obtained as a limit is a triangle $\Delta$ with boundary along
arcs of $\ell_1$, $\cup_k \ell_k$, and $\ell$, with
$\cup_k \ell_k =\{ u= 2k, k\in \Z \}$.
\label{degenerate}
\end{Def}

\begin{Lem}
Suppose given $a_1$ and $\tilde a_1$ in $\M_{Y_1}$, $a \in
\M_{Y,\mu}\backslash j(\M_{Y_1})$, and $a_0\in \M_{Y_0}(\s_k)$. 
Then the holomorphic triangles that appear in the geometric limits of
solutions in the zero-dimensional moduli spaces
$\M^{W_1}_{\ell}(a_1,j(\tilde a_1))$ and 
$\M^{W_0}_k (a,a_0)$ are all degenerate in the limit $\epsilon \to 0$.
\label{tri:deg}
\end{Lem}

\proof
The triangle $\Delta^\epsilon$ has two sides along arcs of the lines
$\{ v-u=1 \}$ and $\{ v=f'(u) \}$ connecting $\vartheta_1$ and $a^-$
and $\vartheta_1$ and $a^+$, respectively, with the parameterization
of Case (e) of Proposition \ref{geom:lim}. If for $\epsilon\to 0$ the
points $\{ \vartheta_1, a^-, a^+ \}$ all lie on the same line $\{
v-u=1 \}$, then the holomorphic map in the limit is a disk with one
side along the arc in $\{ v-u=1 \}$ connecting $a^-$ and $a^+$ and the
other side along an arc in $\partial_\infty(\M_V^*)$ connecting these
same two points. Thus, the limit triangle is degenerate. 

\endproof

Now observe, instead, that when we consider solutions on the cobordism 
that intertwine the generators of the Floer complex for $Y_1$ with
those for $Y_0$, we may have holomorphic triangles that do not
degenerate in the limit when $\epsilon \to 0$.

\begin{Lem}
Suppose given, as before, critical points $a_1$ and $\tilde a_1$ in
$\M_{Y_1}^*$, $a_0(\epsilon) \in 
\M_{Y,\mu}\backslash j(\M_{Y_1})$, and $a_0\in \M_{Y_0}(\s_k)$. 
Consider zero-dimensional moduli spaces of the form
$\M^{W_1}_{\ell}(a_1,a_0(\epsilon))$ and $\M^{W_0}_k (\tilde a_1(\epsilon),
a_0)$, with $\tilde a_1(\epsilon)=j(\tilde a_1)$. 
Then, in general, the limit holomorphic map will still be a triangle. 
\label{tri:nondeg}
\end{Lem}

\proof

Suppose given $a_1$ in $\M_{Y_1}$ and $a_0(\epsilon) \in 
\M_{Y,\mu}\backslash j(\M_{Y_1})$, and consider the geometric limits
of solutions in a zero dimensional moduli space
$\M_{\ell}^{W_1}(a_1,a)$. We still may have holomorphic triangles that
degenerate in the limit $\epsilon \to 0$. By the open mapping theorem, 
this happens whenever two sides of the boundary of the holomorphic
triangle are mapped together as $\epsilon \to 0$. However, since the
pair of points $a_1$
and $a_0 =\lim_\epsilon a_0(\epsilon)$, or $\tilde a_1$ and $a_0$, are
now on two different lines in $H^1(T^2,\R)$, we can also have
non-degenerate holomorphic triangles.  

\endproof

The existence of these non-degenerate holomorphic triangles, that are
not ``thin'' for small  
$\epsilon$, characterizes the difference between the chain maps induced 
by the cobordisms and the homomorphisms of abelian groups defined by
inclusion and projection, as described in Part I \cite{CMW}, as we
shall see when discussing injectivity and surjectivity of the maps in
the exact sequence.

\section{The chain homomorphisms}

In this section we introduce the chain homomorphisms $w_*^1$ and
$w_*^0$.

First we observe that we have an analogue of the formula
\ba
N_{X}(-\s)=(-1)^{(1+b_2^+(X)-b_1(X))/2}N_{X}(\s) 
\label{str:2}
\na 
that holds for compact 4-manifolds. The version we need is given by
the following statement. 

\begin{Lem}
On the manifold $W_1$ we have
$$ N^{W_1}_\ell(a_1,a)= (-1)^{2\, Ind_{\C}(D_{\A})}
N^{W_1}_{-\ell}(a_1,a). $$
Thus, if $\s \in {\cal S}(W_1)$ satisfies
$$ \iota(\s,W_1,a_1,a)= 0, $$
for given $a_1\in \M_{Y_1}^*$ and $a\in \M_{Y,\mu}^*$, then $-\s$ is the 
unique other $\spinc$ structure  in ${\cal S}(W_1)$ which
satisfies $ \iota(-\s,W_1,a_1,a)= 0$. The corresponding invariant
satisfies 
$$ N^{W_1}_{-\s} (a_1,a) = N^{W_1}_{\s} (a_1,a). $$
\label{pos:neg}
\end{Lem}

\proof  The proof for the closed manifold case \cite{Mor}
adapts to this context, with fixed asymptotic values, and compatible
choice of admissible sections. In fact, the sign is given by the
change of orientation, and the orientation is compatible with the
boundary strata.  

\endproof

Notice that the same argument does not extend to the case of the
manifold $W_0$. In fact, in the case of $W_1$ changing $\s$ to $-\s$
does not change the $\spinc$ structure on the boundary $Y\cup Y_1$,
whereas, in the
case of $W_0$, changing $\s$ to $-\s$ amounts to changing the $\spinc$ 
structure on $Y_0$. On the other hand,  on $W_0$ or $W_2$ we simply do
not have solutions  in
$\M_k^{W_0}$ and in $\M_{-k}^{W_0}$ with the same asymptotic values.

Let $a_1=[A_1,\psi_1]$ be a class in ${\cal M}_{Y_1}^*$, and
$a=[A,\psi]$ a class in ${\cal M}_{Y,\mu}^*$, such that the gradings,
assigned according to \cite{CMW}, coincide
$$ \deg_{Y_1}(a_1)= \deg_{Y,\mu} (a). $$
Then we consider  the zero dimensional moduli spaces
${\cal M}^{W_1}_{\s_\ell}(a_1,a)$ with $\s_\ell$, $\ell\geq
0$,  satisfying $\iota(\s_\ell,W_1,a_1,a)=0$. 
Similarly, for $a_0=[A_0,\psi_0]$ in some ${\cal M}_{Y_0}(\s_{k})$
with compatible grading
$$ \deg_{Y,\mu} (a)=\deg_{Y_0,\s_{k}}(a_0), $$
we consider the zero dimensional components
$$ {\cal M}^{W_0}_{\s_k}(a,a_0) \ \ \hbox{with} \ \
\iota(\s_k,W_0,a,a_0)=0. $$  

Recall that we have invariants $N_\ell^{W_1}(a_1,a)$ and
$N^{W_0}_k(a,a_0)$ defined by counting points with the orientation in
${\cal M}^{W_1}_{\s_\ell}(a_1,a)$ and ${\cal M}^{W_0}_{\s_k}(a,a_0)$,
respectively. 

\begin{Def}
We define the map $w_*^1: C_*(Y_1)\to C_*(Y,\mu)$ with matrix elements
$$ \la a , w_*^1(a_1) \ra = N^{W_1}_\ell (a_1,a), $$
with $\ell$ is the unique 
non-negative $\spinc$ structure satisfying 
$\iota(\s_\ell, W_1, a_1,a)=0$.

We define the map $w_*^0: C_*(Y,\mu)\to \oplus_k C_{(*)}(Y_0,\s_k)$
with the matrix coefficients 
$$ \la a_0 , w_*^0(a) \ra =N^{W_0}_k (a,
a_0). $$
\label{w10}
\end{Def}

Thus, we choose to define the map $w_*^1$ using only the ``positive''
$\spinc$ structures, $c_1(\det(\s)) = (2\ell+1) e_1$, with $\ell \geq
0$. 

\begin{Lem}
The maps $w_*^i$ are chain homomorphisms.
\label{chain:homom}
\end{Lem}

\proof  Suppose given a $\spinc$-structure $\s_\ell\in {\cal S}_1$,
with $\ell \geq 0$, satisfying 
$$ \iota(\s_\ell, W_1, a_1,b)=1.$$

We have a compactification of the moduli space
$\M^{W_1}_{\s_\ell}(a_1,b)$ by boundary strata of the form
$$ \bigcup_{a\in \M^*_{Y,\mu}}\M^{W_1}_{\s_\ell}(a_1,a)\times
\hat\M_{Y}(a,b) $$ 
and
$$ \bigcup_{b_1 \in \M^*_{Y_1}} \hat\M_{Y_1}(a_1,b_1) \times
\M^{W_1}_{\s_\ell}(b_1,b), $$ 
as in Proposition \ref{compactif} and Corollary \ref{irredonly}.
Here $a$ and $b_1$ satisfy
$$\deg_{Y,\mu}(a)-\deg_{Y,\mu}(b)=1, $$ 
$$ \deg_{Y_1}(a_1)-\deg_{Y_1}(b_1) =1. $$
The moduli spaces $\hat\M_{Y_1}(a_1,b_1)$ and $\hat\M_{Y}(a,b)$ are
gauge classes of flow lines on $Y_1$ and $Y$ respectively, modulo the action
of $\R$ by translation.
The counting of boundary points of the 1-dimensional
$\M^{W_1}_{\s_\ell}(a_1,b)$, with the orientation, gives the relation
$$ \sum_{a} N^{W_1}_{\s_\ell}(a_1,a) n_{Y}(a,b) - \sum_{b_1}
n_{Y_1}(a_1,b_1) N^{W_1}_{\s_\ell}(b_1,b) =0. $$
Notice that $\M^{W_1}_{\s_\ell}(a_1,a)$ and
$\M^{W_1}_{\s_\ell}(b_1,b)$ are zero-dimensional, thus we have 
$$ \iota(\s_\ell,W_1,a_1,a)=\iota(\s_\ell,W_1,b_1,b)=0. $$
Thus, counting points in these moduli spaces gives exactly the
counting of the matrix elements of $w_*^1$. 
This proves the relation
$$ w_*^1\circ \partial_{Y_1} =\partial_{Y,\mu} \circ w_*^1. $$

The result for $W_0$ is analogous, except for the fact that some care
is needed in the case of the trivial $\spinc$-structure $\s_0$. 
Let us first consider the case of non-trivial $\s_k$ first. Suppose given
a 1-dimensional $\M^{W_0}_{\s_k}(a,b_0)$.
In particular, we have 
$$ \deg_{Y,\mu}(a)= \deg_{Y_0}(b_0) +1, $$
where $\deg_{Y_0}$ is the $\Z$-lift of the $\Z_{2k}$-relative
grading on $\M_{Y_0}(\s_k)$, induced by the grading on $\M_{Y,\mu}$,
cf.\cite{CMW}, \cite{MW4}. 

The boundary strata in the compactification of
$\M^{W_0}_{\s_k}(a,b_0)$ are given by
\ba \label{comp:W0} \begin{array}{l} \bigcup_{a_0\in \M_{Y_0}(\s_k)}
\M^{W_0}_{\s_k}(a,a_0)\times 
\hat\M_{Y_0,\s_k}^{(0)}(a_0,b_0)\\[3mm]
\bigcup_{b\in \M_{Y,\mu}^*} \hat\M_{Y,\mu}(a,b)\times
\M^{W_0}_{\s_k}(b,b_0). \end{array} \na
Here $b$ and $a_0$ satisfy 
$$ \deg_{Y,\mu}(a)-\deg_{Y,\mu}(b)=1, $$
\ba \deg_{Y_0}(a_0)- \deg_{Y_0}(b_0) =1, \label{Z=}\na
in the integer grading.
Notice that, as in \cite{MW4}, (\ref{Z=}) can be derived by
observing that, in the geometric limits for $r\to\infty$, the
flowlines that contribute to the compactification define a
contractible path $[A''(t),0]$ in $\M_{\nu(K)}\subset
\chi_0(T^2,\nu(K))$. 
Thus, with the notation of \cite{MW4}, we pick only the minimal energy
component  
$\M_{Y_0,\s_k}^{(0)}(a_0,b_0)$ in the moduli space
$\M_{Y_0,\s_k}(a_0,b_0)$. 
This component of the moduli space is
exactly the one which defines the boundary operator on the $\Z$-lift
$C_{(*)}(Y_0,\s_k)$  of the $\Z_{2k}$ graded complex
$C_{*}(Y_0,\s_k)$, as analyzed in \cite{MW4}.

The argument used here, based on the geometric limits of solutions,
can be generalized to the case of the trivial $\spinc$-structure $\s_0$.
Again, suppose given $\M^{W_0}_{\s_0}(a,b_0)$, with 
$$ \iota(\s_0,W_0,a,b_0)=1. $$

In this case we want to show that the strata
$$ \M^{W_0}_{\s_0}(a,a_0)\times \hat\M_{Y_0,\s_0}(a_0,b_0) $$
which occur in the compactification, with 
$$ \deg_{Y_0,\s_0}(a_0)-\deg_{Y_0,\s_0}(b_0)=1, $$
only contain flow lines in the component of minimal energy of
$\M_{Y_0,\s_0}(a_0,b_0)$. Again, the result follows from the fact that the
path $[A''(t),0]$ in the geometric limits is contractible, for any
flow-line in $\M_{Y_0,\s_0}(a_0,b_0)$ that arises in the
compactification of $\M^{W_0}_{\s_0}(a,b_0)$.
This precisely identifies the right component of the boundary operator in the
Floer complex for $(Y_0, \s_0)$.

Thus, (\ref{comp:W0}) implies that the
map 
$$ w_*^0: CF_*(Y,\mu)\to \oplus_k CF_{(*)}(Y_0,\s_k)$$
is a chain homomorphism precisely when the complex
$CF_{(*)}(Y_0,\s_k)$ is endowed with the $\Z$-grading described in
\cite{MW4} and, in the case of $\s_0$, the Floer homology of $Y_0$ is
defined as 
$$ HF_{(*)}(Y_0,\s_0, \Z[[t]])|_{t=0} $$ on
$(Y_0,\s_0)$, as described in \cite{MW4}. 

\endproof

\subsection{The composite map}

Let $N^W_{\ell, k}(a_1,a_0)$ denote the invariant obtained by counting
solutions in moduli space $\M^W_{\ell, k}(a_1,a_0)$, for the unique
choice of the $\spinc$ structure $\s_{\ell,k}$ such that
$\M^W_{\ell, k}(a_1,a_0)$ is zero dimensional. 

\begin{Lem}
Suppose given $a_1$ in $\M_{Y_1}^*$ and $a_0$ in
$\M_{Y_0}(\s_k)$. Let $\ell \geq 0$ be the unique non-negative
$\spinc$-structure satisfying $$\iota(\s_{\ell,k},W,a_1, a_0)=0.$$ 
Then the composite map $w^0_* \circ w^1_*$ is given by
$$ \la w^0_* \circ w^1_*(a_1), a_0 \ra = N^W_{\ell,k}(a_1, a_0). $$
\label{split:W}
\end{Lem}

\proof

The composite map $w^0_* \circ w^1_*$ has matrix  elements
$$  \sum_{a\in \Theta}  N^{W_1}_{\ell}(a_1, a) N^{W_0}_k(a, a_0), $$
where the sum is over the set 
$$ \begin{array}{ll}\Theta=& \{ a\in {\cal M}_{Y,\mu}^* |
\iota(\s_{\ell},W_1,a_1,a)=0, 
\iota(\s_k,W_0,a,a_0)=0, \\ & \deg_{Y_1}(a_1)=\deg_{Y,\mu}(a)=\deg_{Y_0}(a_0) 
\}. \end{array} $$ 
On the other hand, we have
$$ N^W_{\ell,k}(a_1, a_0)=\sum_{a\in \Theta_1} N^{W_1}_{\ell}(a_1, a)
N^{W_0}_k(a, a_0), $$
where now the sum is over the set
$$ \Theta_1=\{ a\in {\cal M}_{Y,\mu}^* | \iota(\s_{\ell},W_1,a_1,a)=0,
\iota(\s_k,W_0,a,a_0)=0 \}. $$

Notice that, given any $\s_{\ell,k}$ on $W$ which satisfies 
$$ \iota( \s_{\ell,k},W,a_1,a_0)=0, $$
and such that we have
$$ \M^W_{\s_{\ell,k}}(a_1,a_0)\neq \emptyset, $$
there exists some $a\in \M_{Y,\mu}^*$ such that we have
$$\iota(\s_{\ell},W_1,a_1,a)=0 \ \ \hbox{and} \ \ 
\iota(\s_k,W_0,a,a_0)=0 . $$ This follows by
stretching the cylinder $Y\times [-T,T]$ in the composite
cobordism: the condition
$$ \M^W_{\s_{\ell,k}}(a_1,a_0)\neq \emptyset $$
ensures the existence of a limiting translation invariant solution on
$Y\times \R$. The argument is similar to the one used in \cite{KM}.

Thus, we only need to prove that we have
$$ \sum_{a\in  \Theta_1\backslash \Theta} N^{W_1}_{\ell}(a_1, a)
N^{W_0}_k(a, a_0) =0. $$
Suppose not. Then we have an element $a\in \Theta_1\backslash \Theta$
such that the moduli spaces
$$ \M^{W_1}_{\ell}(a_1, a)  \ \ \hbox{ and } \ \ \M^{W_0}_k(a, a_0) $$
are zero-dimensional and non-empty. The point $a$ will be in
$j(\M_{Y_1}^*)$, or in $\M_{Y,\mu}^*\backslash j(\M_{Y_1}^*)$. We
consider the first case. The proof in the second case is completely
analogous. We have $a=j(\tilde a_1)$ and $\deg_{Y_1}(a_1)\neq
\deg_{Y_1}(\tilde a_1)$. Consider the geometric limits of the
solutions in $\M^{W_1}_{\ell}(a_1,j(\tilde a_1))$. As we discuss in
greater detail in Lemma \ref{constant} and Lemma \ref{rel:2}, as we
let $\epsilon\to 0$, the geometric limits of solutions in
$\M^{W_1}_{\ell}(a_1,j(\tilde a_1))$ define geometric limits of
solutions in $\M_{Y_1}(a_1,\tilde a_1)$. In particular, if the moduli
space $\M^{W_1}_{\ell}(a_1,j(\tilde a_1))$ is non-empty and
zero-dimensional, also the moduli space $\M_{Y_1}(a_1,\tilde a_1)$
will be non-empty and zero-dimensional, which contradicts the
assumption that $\deg_{Y_1}(a_1)\neq \deg_{Y_1}(\tilde a_1)$.

This implies that we have an identification of the
zero-dimensional components
\ba \M^W_{\s_{\ell,k}}(a_1,a_0) \cong  \cup_{a \in \Theta }
\M^{W_1}_{\s_{\ell}}(a_1,a) \times \M^{W_0}_{\s_k}(a,a_0), \label{0-comp}\na
for large $T\geq T_0$ in the cylinder $Y\times [-(T-T_0),(T-T_0)]$ in the
composite cobordism $W$. This follows from the gluing theorem
of \cite{MW}, cf. \cite{CW} \cite{MW4}, together with the previous argument.

\endproof

\section{Injectivity and Surjectivity}

In the previous section we have constructed a sequence 
$$ 0\to C_*(Y_1)\stackrel{w^1_*}{\to}
C_*(Y,\mu)\stackrel{w^0_*}{\to} \oplus_k C_{(*)}(Y_0,\s_k) \to 0. $$
We now proceed to show exactness in the first and last place, namely
injectivity of $w^1_*$ and surjectivity of $w^0_*$.

We want to give a better description of the components
$$ \la j(a_1'), w_*^1(a_1) \ra, $$
with $a_1$ and $a_1'$ in $\M_{Y_1}^*$,
and
$$ \la a_0, w_*^0( a)\ra, $$
$a\in \M_{Y,\mu}^*\backslash j(\M_{Y_1}^*)$, and $a_0 \in \cup_k
\M_{Y_0}(\s_k)$. 
For simplicity, let us introduce the following notation.
Let $\{ a^{(1)}_i \}_{i=1,\ldots m}$ be the elements in ${\cal
M}_{Y_1}^*$ and let $\{ a^{(0)}_j \}_{i=m+1,\ldots n}$ be the elements in
$\cup_{\s_{i_k}}{\cal M}_{Y_0}(\s_{i_k})$.
Then the $n$ elements in ${\cal M}_{Y,\mu}^*$ can be identified with the
union of these two sets of points. More precisely, if $\epsilon >0$ is
the parameter used in the 
definition of the function $f'$ in the construction of the surgery
perturbation, then the moduli space $ {\cal M}_{Y,\mu}^*$ can be
identified with a collection of points
$$ {\cal M}_{Y,\mu}^* = \{ a^{(1)}_i(\epsilon) \}_{i=1,\ldots m} \cup \{
a^{(0)}_j(\epsilon) \}_{j=m+1,\ldots n}. $$
For $\epsilon$ small enough, there is a bijection
\ba
{\cal M}_{Y,\mu}^* \cong \{ a^{(1)}_i \}_{i=1,\ldots m} \cup \{ a^{(0)}_j
\}_{j=m+1,\ldots n},  \label{MY} \na
which is compatible with the grading \cite{CMW}.

In the following we shall use the notation
(\ref{MY}) which identifies the elements in ${\cal M}_{Y,\mu}$ with
elements in the other two moduli spaces.Whenever it is crucial to
distinguish between these moduli spaces, we shall use the notation
$a_1 \in \M_{Y_1}^*$, $a_0\in 
\M_{Y_0}(\s_k)$, and $a\in \M_{Y,\mu}^*$, with $a=j(a_1)$ or
$\pi(a)=a_0$.  We hope this will not cause any confusion.

We want to describe solutions in the moduli spaces
$$ \M_\ell^{W_1} (a^{(1)}_i, a^{(1)}_q) $$
and
$$ \M_k^{W_0} (a^{(0)}_j,a^{(0)}_p ), $$
by describing these moduli spaces as a gluing of
solutions on the trivial cobordism $V\times \R$ and  solutions on the
regions  $\hat W_i(\nu(K))$. The techniques involved in the splitting
and gluing of solutions on the cobordisms are analogous to the ones
developed in \cite{MW2} to analyze the splitting and gluing of flow
lines on the trivial cobordism.

We are assuming here that the asymptotic values also satisfy the
condition 
\ba \label{same:deg} \begin{array}{c} 
\deg_{Y_1}(a^{(1)}_i)=\deg_{Y,\mu}( a^{(1)}_q), \\
\deg_{Y,\mu}( a^{(0)}_j )=\deg_{Y_0,\s_k} (a^{(0)}_p).
\end{array} \na

In particular, we are going to prove the fundamental relations
$$ \la a^{(1)}_i, w_*^1(a^{(1)}_q \ra =\delta_{iq} $$
for all $a^{(1)}_i$ and $ a^{(1)}_q$ in $\M_{Y_1}^*$ and 
$$ \la a^{(0)}_j, w_*^0(a^{(0)}_p) \ra =\delta_{jp} $$
for all $a^{(0)}_j$ and $ a^{(0)}_p$ in $\cup_k \M_{Y_0}(\s_k)$.

In the following we shall discuss the case of the manifold $W_1$. The
case of $W_0$ is analogous.
Consider the manifold 
$$ \hat W_1(r)= (V_r \times \R) \cup {\cal V}(r) \cup \hat {\cal V}(r), $$
as in (\ref{W:parts}).

We use the description of the geometric limits of solutions on $W_1$
given in Proposition \ref{geom:lim}, together with the analysis of
\cite{MW2} of the geometric limits of flow lines, 
in order to analyst solutions in $\M^{W_1}_\ell (a^{(1)}_i,a^{(1)}_i)$
and in $\M^{W_0}_k (a^{(0)}_j ,a^{(0)}_j )$. 

Suppose given $a_1$, an element of $\M_{Y_1}^*$. For large $r\geq r_0$
we represent $a_1$ as in \cite{CMW},
$$ a_1(r)=[(\tilde A',\tilde \psi')\#_{\tilde a_\infty''}^r (\tilde
a_\infty'',0)]. $$

\begin{Lem}
There is a unique solution $[\A_r,\Psi_r]$  in the zero-dimensional
moduli space 
$$\M^{W_1(r)}_\ell (a_1 (r),j(a_1(r))),$$
for large $r \geq r_0$ and for small enough $\epsilon >0$, where
$\epsilon$ is the parameter used in the definition of the surgery
perturbation. 
\label{constant}
\end{Lem}

\proof
We write $a_1(r)=[A(r),\psi(r)]$, with
$$ [A(r),\psi(r)]=[(\tilde A',\tilde \psi')\#_{\tilde a_\infty''}^r
(\tilde a_\infty'',0)] $$ 
on $Y_1(r)$, for large $r\geq r_0$, and
$$ j(a_1(r))=j[A(r),\psi(r)]=[(A',\psi')\#_{a_\infty''}^r
(a_\infty'',0)]. $$  

In the limit $\epsilon \to 0$, the asymptotic values $a_\infty''$ and
$\tilde a_\infty''$ coincide, and the elements $[A',\psi']$ and
$[\tilde A',\tilde \psi']$ in $\M_V^*$ also coincide, because of the
result of Lemma \ref{tri:deg}, which shows that the holomorphic
triangle degenerates. Thus, when
$\epsilon \to 0$ there is a unique finite energy solution on $V\times
\R$ with the same asymptotic value $[A',\psi']$ at $t\to\pm\infty$,
given by the constant flow $[A',\psi']$.

In fact, if we had a non-constant finite energy solution
$(\A',\Psi')$ on $V\times \R$ with the same limits at $t\to\pm\infty$,
then, using the gluing theorem of \cite{MW2}, we could glue this
$(\A',\Psi')$ 
along the flat connection $a_\infty''$ to a solution on $Y_1\times \R$,
extending it as the reducible solution $(a_\infty'',0)$ on $\nu(K)\times
\R$. This solution would be a non-trivial flow line
connecting the critical point $a_1$ to itself in the configuration
space on $Y_1$, but the moduli space $\hat \M_{Y_1}(a_1,a_1)$ 
is generically empty for dimensional reasons.

Thus, in the limit $\epsilon \to 0$, there is a unique solution
$[\A_r,\Psi_r]$ on $W_1(r)$ with the same asymptotic value at the two
ends, obtained by gluing along
the asymptotic value $a_\infty''$ the constant flow $[A',\psi']$ with
the unique perturbed ASD equation on $W_1(\nu(K))$ which extends
$a_\infty$. 

Since the moduli space $$\M^{W_1(r)}_\ell (a_1 (r),j(a_1(r)))$$ is
discrete for all $0< \epsilon \leq \epsilon_0$, we obtain that, for
$\epsilon$ small enough, there is a unique solution.

\endproof

\begin{Cor}
We have 
$$ \la a^{(1)}_i, w_*^1(a^{(1)}_i) \ra =1 $$
for all $a^{(1)}_i$ in $\M_{Y_1}^*$, and the analogous
$$ \la a^{(0)}_j, w_*^0(a^{(0)}_j) \ra =1 $$
for all $a^{(0)}_j$ in $\cup_k \M_{Y_0}(\s_k)$.
\label{fund:1}
\end{Cor}

We also have the following result.

\begin{Lem}
Consider the moduli space
$\M^{W_1}_\ell(\tilde a_1 , j(a_1))$, for large enough $r\geq r_0$
and with two critical points $a_1 \neq \tilde a_1$ in $\M_{Y_1}^*$.
Assume as before that $0< \epsilon \leq \epsilon_0$ is the parameter
used in the construction of the surgery perturbation.
Then, for $\epsilon$ small enough, we have 
$$ \M^{W_1}_\ell(\tilde a_1 , j(a_1)) = \emptyset $$
if $\deg_{Y_1}(a_1)= \deg_{Y_1}(\tilde a_1)$.
\label{rel:2}
\end{Lem}

\proof
Suppose given a solution $[\A_r,\Psi_r]$ in
$$ \M^{W_1}_\ell(\tilde a_1 , j(a_1)),$$
with $a_1 \neq \tilde a_1$ in $\M_{Y_1}^*$ and
$$ \deg_{Y_1}(a_1)= \deg_{Y_1}(\tilde a_1). $$
According to the result of Proposition \ref{geom:lim}, we have
geometric limits as $r\to\infty$. We show that we can construct from
these geometric limits a non-trivial flow line in the
configuration space over $Y_1$ which connects the two points
$a_1$ and $\tilde a_1$ in $\M_{Y_1}^*$. This will contradict the assumption
that 
$$ \deg_{Y_1}(a_1)= \deg_{Y_1}(\tilde a_1). $$
In fact, this assumption implies that the moduli space of flow lines is
generically empty for dimensional reasons.
 
We can write the endpoints $a_1$ and $\tilde a_1$ as
$$ j(a_1)= [(A ',\psi ')\#_{a_\infty}^r (a_\infty, 0)], $$
and
$$ \tilde a_1= [(\tilde A ',\tilde \psi ')\#_{\tilde a_\infty}^r
(\tilde a_\infty, 0)]. $$
Here $a_\infty$ is a flat connection on $T^2$ satisfying the relation
$v=f'(u)$, with $f'$ depending on the choice of a small $\epsilon >0$,
and $\tilde a_\infty$ is a flat connection on $T^2$ satisfying the
relation $v-u=1$.

Suppose that we have $\M^{W_1}_{\ell}(\tilde a_1,j(a_1))\neq
\emptyset$. Then,
in the limit $\epsilon \to 0$, the geometric limits on $\hat W_1(\nu(K))$
consist of the cases (c), (d), and (e) of Proposition \ref{geom:lim},
where the connections $a_\infty''$ and $\tilde a_\infty ''$
coincide, because of the result of Lemma \ref{tri:deg} on the
degenerate holomorphic triangles. Therefore, these geometric limits
determine a preglued solution on 
$\nu(K)\times \R$ inside $Y_1\times \R$ which can be glued to the
geometric limits on $V\times \R$ to determine a flow line on
$Y_1(r)\times \R$. Since we know that $\M_{Y_1(r)\times \R}(\tilde
a_1,a_1)=\emptyset$ if $\deg_{Y_1}(a_1)= \deg_{Y_1}(\tilde a_1)$,
we have shown that, for $\epsilon \to 0$, we also have 
$$ \M^{W_1}_\ell(\tilde a_1 , j(a_1)) =\emptyset, $$
hence the relation
$$ \langle w_*^1 (\tilde a_1), j(a_1) \rangle =0  $$
holds for small enough $\epsilon$, if $a_1 \neq \tilde a_1$ in
$\M_{Y_1}^*$. 

\endproof

We derive from Lemma \ref{rel:2} the second fundamental relation.

\begin{Cor}
We have
$$ N^{W_1}_\ell(a^{(1)}_i, a^{(1)}_q) =0 $$
whenever $i\neq q$ and 
$$ \deg_{Y_1}(a^{(1)}_i)= \deg_{Y_1}(a^{(1)}_q). $$
Similarly, for $W_0$ we have
$$ N^{W_0}_k (a^{(0)}_j, a^{(0)}_p) =0 $$
whenever $j\neq p$ and 
$$ \deg_{Y_0}(a^{(0)}_j)= \deg_{Y_0}(a^{(0)}_p). $$
\label{fund:2}
\end{Cor}

Notice how the previous results do not give any information about the
components of the maps $w_*^1$ and $w_*^0$ that interchange the
critical points $a^{(1)}_i$ and $a^{(0)}_j$. In fact,  in
this case there are in general non-degenerate triangles, as 
$\epsilon \to 0$ in the surgery perturbation, cf. the result of Lemma
\ref{tri:nondeg}. In other words, the presence of these non-degenerate
holomorphic triangles measures the difference between the chain maps
$w_*^1$ and $w_*^0$ and the group homomorphisms $j$ and $\pi$ defined
by the inclusion and projection, in the identification of the
generators on $Y$ with generators on $Y_1$ and $Y_0$.

\section{Exactness in the middle term}

Recall that the results of \cite{CMW} imply that the ranks of the Floer
complexes in the sequence
$$ 0\to C_*(Y_1) \to C_*(Y,\mu) \to \oplus_k C_{(*)}(Y_0,\s_k) \to 0
$$ 
are as prescribed for the existence of an exact sequence. Moreover, in 
the previous section we have proved injectivity of the first map and
surjectivity of the last. Now we analyze the middle term. 

The result of Lemma \ref{split:W}, together with Corollary
\ref{fund:1} and Corollary \ref{fund:2}, yields the following.

\begin{Lem}
Suppose given $a\in \M_{Y,\mu}\backslash j(\M_{Y_1})$ and $a_1\in
\M_{Y_1}$ and let $\pi$ be the identification $\pi: \M_{Y,\mu}\backslash
j(\M_{Y_1}) \to \cup_k \M_{Y_0}(\s_k)$.
The coefficients of the composite map $w_*^0 \circ w_*^1$ satisfy the
relation 
\ba N^W_{\ell,k}(a_1,\pi(a)) = N^{W_1}_{\ell}(a_1,a) + N^{W_0}_k
(j(a_1),\pi(a)). \label{count:W} \na
\label{comp:rel}
\end{Lem}

The main purpose of this section is to show that the counting in
(\ref{count:W}) is zero. We then verify that this is sufficient to
prove exactness in the middle term.
Establishing the relation $w_*^0\circ w_*^1=0$ depends again
essentially on the analysis of the geometric limits of solutions 
on the cobordisms, following the technique of
\cite{MW2}. We need a preliminary discussion of the geometric limits
on $V\times \R$ which completes the results of Part II, \cite{MW2}.

\subsection{The moduli space on $V\times \R$}

Consider elements $a_i^{(1)}\in \M_{Y_1}^*$ and $a_j^{(0)}(\epsilon)
\in \M_{Y,\mu}^*$, which we can write as
$$ a_i^{(1)}=[(A_i^-,\psi_i^-)\#(a_{\infty,i}^-,0)] $$
$$ a_j^{(0)}(\epsilon)= [(A_j^+(\epsilon),\psi_j^+(\epsilon))\#
(a_{\infty,j}^+,0)]. $$

Assume that we have solutions $[\A_1(r), \Psi_1(r)]$ in  
$\M_\ell^{W_1(r)}\bigl(a_i^{(1)}, a_j^{(0)}(\epsilon)\bigr)$, for all
sufficiently large $r \geq r_0$. Then these solutions define geometric 
limits as in Proposition \ref{geom:lim} (a) - (h).

In particular, we list here separately the limits on $V\times \R$. Our
purpose now is to simplify and group together in a more efficient way
the information on the geometric limits given in Proposition
\ref{geom:lim}.

\begin{Rem}
\label{geom:lim:V}
A family of solutions $[\A_1(r), \Psi_1(r)]$ in  
$\M_\ell^{W_1(r)}\bigl(a_i^{(1)}, a_j^{(0)}(\epsilon)\bigr)$ defines
the following limits on $V\times \R$:
 
\noindent 
 (a). A finite energy solution $[\A',\Psi']^\epsilon$ of the perturbed
equations (\ref{4SWP}) on $V \times \R$, with a radial limit
$a_\infty(\epsilon)$ in $\partial_\infty (\M_V^*)\subset
\chi_0(T^2,V)$, and with temporal limits $[A,\psi]_1^\epsilon$ and
$[\tilde A, \tilde \psi]_1^\epsilon$ in
$$ \partial_\infty^{-1}(a_\infty(\epsilon)) \subset \M_V^*. $$

\noindent 
(b). Two paths $[A(t),\psi(t)]_1^\epsilon$ in $\M_V^*$, for
$t\in [-1,0)$ and $t\in (0,1]$, with
$$ [A(-1),\psi(-1)]_1^\epsilon = [A_i^-, \psi_i^-] $$
$$ \lim_{t\to 0-} [A(t),\psi(t)]_1^\epsilon =
[A,\psi]_1^\epsilon $$
$$ \lim_{t\to 0+} [A(t),\psi(t)]_1^\epsilon = [\tilde A,
\tilde \psi]_1^\epsilon $$
$$ [A(1),\psi(1)]_1^\epsilon = [ A_j^+(\epsilon) ,
\psi_j^+(\epsilon)]. $$

These paths induce a continuous, piecewise smooth path 
$$ a_1^\epsilon(t) \subset \partial_\infty (\M_V^*)\subset
\chi_0(T^2,V)$$  satisfying
$$ a_1^\epsilon(t)=\partial_\infty
[A(t),\psi(t)]_1^\epsilon, $$
with
$$ a_1^\epsilon (-1)=a_{\infty,i}^-  \ \  a_1^\epsilon
(0)=a_\infty(\epsilon) \ \  a_1^\epsilon (1)=
a_{\infty,j}^+(\epsilon). $$ 

As $\epsilon\to 0$, these geometric limits define paths
$[A(t),\psi(t)]$ and 
$a(t)$ with similar properties, and with
$$ a(-1)=a_{\infty,i}^- \ \  a(0)=a_\infty \ \ 
a(1)=a_{\infty,j}^+, $$
in $\partial_\infty (\M_V^*)\subset
\chi_0(T^2,V)$, with $a_\infty=\lim_\epsilon a_\infty(\epsilon)$.

(c) Moreover, we have a holomorphic triangle in $H^1(T^2,\R)$ with
vertices 
$$\{ a^-_{\infty,i}, \vartheta_1, a^+_{\infty,j}(\epsilon) \}$$
and sides given by parameterized arcs along the lines $\ell_1=\{ v-u=1
\}$, $\ell_\mu=\{ v =f'(u) \}$ and by
$$ \{ a_1^\epsilon(t) \}\subset \ell=\partial_\infty(\M_V^*). $$
This is obtained from the non-uniform limits of Proposition
\ref{geom:lim}. 

\end{Rem}

We now describe how to assemble together these geometric limits in a
suitable moduli space. This will be useful in the following subsection,
in the proof of Theorem \ref{W_1:W_0} that establishes the exactness
in the middle term.

Let $a_\infty$ be an element in $\chi_0(T^2,V)$.
We define the configuration space 
$$ {\cal A}_{k,\delta}(V\times \R, a_\infty) $$
as follows. 

We can write a pair $(\A,\Psi)$ of a $U(1)$-connection and a spinor on
$V\times \R$ in the form 
$$ \A= a(w,s,t) + f(w,s,t) ds + h(w,s,t) dt $$
in the region $T^2\times [0,\infty)\times \R$ inside $V\times \R$.  
A pair $(\A,\Psi)$ as above is in the configuration space
${\cal A}_{k,\delta}(V\times \R, a_\infty)$ if $(\A,\Psi)$ is in $L^2_k$ on
$V\times \R$. Moreover, we also require that, after the change of
coordinates $s+it=e^{\rho +i\theta}$, and the corresponding change of
variables 
\ba\begin{array}{l} a(w,\rho,\theta)=a(w,e^{\rho +i\theta}) \\[2mm]
 f(w,\rho,\theta)=e^{-\rho} \cos\theta \ h(w,e^{\rho
+i\theta})-e^{-\rho}\sin\theta \ f(w,e^{\rho +i\theta})\\[2mm] 
 h(w,\rho,\theta)=e^{-\rho}\cos\theta \ f(w,e^{\rho+i\theta})
+e^{-\rho}\sin\theta \ h(w,e^{\rho+i\theta}),
\end{array}\label{change:var} \na
we have
$$ (a-a_\infty,h,f-f_0,\alpha,\beta) \in L^2_{k,\delta}(T^2\times [-\pi/2,
\pi/2] \times [\rho_0,\infty)), $$
where $f_0$ is a constant. The $L^2_{k,\delta}$ norm we consider
is defined by
$$ \| F \|_{L^2_{k,\delta}}= \| e_\delta \cdot F \|_{L^2_k}, $$
where $e_\delta$ is a smooth non-negative function satisfying
$$ e_\delta(w,\rho,\theta)=\exp(\delta e^\rho), $$
for $(w,\rho,\theta)$ in the range 
$$T^2\times [\rho_0,\infty)\times
[-\pi/2,\pi/2].$$ 

We have a group of gauge transformations acting on this configuration
space, namely the group
$$ {\cal G}_{V\times \R, k+2,\delta} $$
given by maps $\lambda : V\times \R \to U(1)$ of the form
$\lambda=e^{i\ell}$ with $\ell: V\times \R \to \R$, such that
$\ell-\ell_\infty$ is in $L^2_{k+2,\delta}$ on the domain
\ba T^2\times [\rho_0,\infty)\times
[-\pi/2,\pi/2] \subset V\times \R.\label{asympt:region} \na 
Here $\lambda_\infty=e^{i\ell_\infty}$ is a gauge transformation on
$T^2$ that extends to $V$, that is,  $\lambda_\infty\in {\cal G}_V$.

Recall that, by the analysis of \cite{MW2}, any finite energy solution
$(\A,\Psi)$ on $V\times \R$ has the property that 
in radial gauge, $(\A,\Psi)$ has limits 
$$ a_\infty(\theta)= \lambda(\theta) a_\infty, $$
with $\theta\in [-\pi/2, \pi/2]$ and
with $\lambda(\theta)$ a family of gauge transformations on $T^2$ that 
extend to gauge transformations on $V$, and 
$[ a_\infty ]$ a fixed gauge class of flat connections on $T^2$.
With a slight abuse of notation, we write $a_\infty$ for this
gauge class. Thus, we can represent all finite energy solutions by
elements in some configuration space ${\cal A}_{k,\delta}(V\times \R,
a_\infty)$. 

Given a fixed element $\Gamma_{a_\infty}=(\A,\Psi)$ in the
configuration space ${\cal A}_{k,\delta}(V\times \R, a_\infty)$, we
define the slice at $\Gamma_{a_\infty}$ as
$$ {\cal S}_{\Gamma_{a_\infty}}=\{ (\tilde a, \tilde \Phi)\in {\cal
A}_{k,\delta}(V\times \R, a_\infty) | 
G^*_{\Gamma_{a_\infty}}(\tilde a, \tilde \Phi)=0 \},   $$
where the operator $G^*_{\Gamma_{a_\infty}}$ is the $L^2_{\delta}$
adjoint of the infinitesimal gauge action $G_{\Gamma_{a_\infty}}$.
In the region (\ref{asympt:region}) the elements $(\tilde a, \tilde
\Phi)$ can be written in the form
$$ (\tilde a, \tilde \Phi)=(a,h,f,\alpha,\beta), $$
with
$$ \tilde a= a+ f\, ds +h\, dt \ \hbox{ and } \ \tilde
\Phi=(\alpha,\beta), $$
as above.

We define the moduli space $\M_{V\times \R}(a_\infty)$ as the set of
solutions of the Seiberg-Witten equations in ${\cal
A}_{k,\delta}(V\times \R, a_\infty)$, modulo the action of the gauge
group ${\cal G}_{V\times \R, k+2,\delta}$. We denote by
$\hat\M_{V\times \R}(a_\infty)$ the {\em balanced energy} moduli space,
namely the elements in $\M_{V\times \R}(a_\infty)$ satisfying 
$$ \begin{array}{l} 
\int_{-\infty}^0 (\| \partial_t A(t) \|_{L^2(V)}^2 + \| \partial_t
\psi(t) \|^2_{L^2(V)})dt = \\[2mm]
\int_0^{+\infty}(\| \partial_t A(t) \|_{L^2(V)}^2 + \| \partial_t
\psi(t) \|^2_{L^2(V)})dt, \end{array} $$
for a temporal gauge representative $(A(t),\psi(t))$.

The linearization ${\cal L}_{(\A,\Psi)}$ at a solution $(\A,\Psi)$ in
the slice  
${\cal S}_{\Gamma_{a_\infty}}$ is given by the operator
\ba {\cal L}_{(\A,\Psi)}(\tilde a, \tilde \Phi)= \left\{
\begin{array}{l} d^+ \tilde a -\frac{1}{2} Im (\Psi\cdot \tilde \Phi) \\[2mm]
D_{\A} \tilde \Phi + \tilde a. \Psi \\[2mm]
G^*_{\Gamma_{a_\infty}}(\tilde a, \tilde \Phi) \end{array}
\right. \label{lin:slice:VR} \na

The virtual dimension of the moduli space $\hat \M_{V\times
\R}(a_\infty)$ at a solution $(\A,\Psi)$ in the slice ${\cal
S}_{\Gamma_{a_\infty}}$ is given by 
$$ Index({\cal L}_{(\A,\Psi)}). $$

Assuming that the element $a_\infty \neq \vartheta$ is away from the bad
point $\vartheta$ in the character variety of $T^2$, we know, by the result of
Section 3.2 of \cite{MW2}, that all finite energy solutions $(\A,\Psi)$
have a uniform exponential decay in radial gauge, hence they can be
regarded as elements of the moduli spaces $\hat \M_{V\times
\R}(a_\infty)$ introduced here, for some $\delta$ which depends only
on $a_\infty$. 

Recall also that in Proposition 3.4 of Part II we proved that, given a
finite energy solution $(\A,\Psi)$ in ${\cal A}_{k,\delta}(V\times
\R,a_\infty)$ of the Seiberg-Witten equations on
$V\times \R$, with asymptotic value in the gauge class of
$a_\infty$, by applying a gauge transformation $\lambda(\theta)$ in
${\cal G}_{V\times \R, k+2,\delta}$, we obtain a solution $\lambda
(\A,\Psi)$ with $h\equiv 0$, and $f-f_0\equiv 0$. In particular, when
inverting the change of variables (\ref{change:var}), the condition
that $\lambda (\A,\Psi)$ is in $L^2_k$ gives $f_0\equiv 0$. In
particular the resulting solution is in a temporal gauge, hence we
can write $\lambda (\A,\Psi)$ as $(A(t),\psi(t))$.
Thus, given a solution $(\A,\Psi)$ in ${\cal A}_{k,\delta}(V\times
\R,a_\infty)$, up to gauge transformations, we obtain two classes
$[A,\psi]$ and $[\tilde A,\tilde \psi]$ in 
$$ \partial_\infty^{-1}(a_\infty) \subset \M^*_V $$  
defined by the asymptotic values as $t\to\pm\infty$ of the temporal
gauge representative $(A(t),\psi(t))=\lambda (\A,\Psi)$, cf. the
result of Lemma 3.9 of \cite{MW2}. 

Thus, we can break the moduli space $\hat \M_{V\times \R}(a_\infty)$
into a union of components  
$$ \hat \M_{V\times \R}(a_\infty) = \bigcup_{[A,\psi],[\tilde A,\tilde
\psi]}\hat \M_{V\times \R}( [A,\psi],[\tilde A,\tilde
\psi], a_\infty), $$
with 
$$ [A,\psi],[\tilde A,\tilde\psi] \in \partial_\infty^{-1}(a_\infty)
\subset \M^*_V. $$  
We can rephrase the calculation of the virtual dimension as follows.

\begin{Pro}
Let $[\A,\Psi]$ be a gauge class in $$\hat \M_{V\times \R}(a_\infty).$$ Let
$(A(t),\psi(t))$ be a temporal gauge representative of $[\A,\Psi]$,
which satisfies
$$ \lim_{t\to -\infty} [A(t),\psi(t)]=[A,\psi], $$
$$ \lim_{t\to \infty} [A(t),\psi(t)]=[\tilde A,\tilde \psi]. $$
Let $H_{A(t),\psi(t)}$ be the operator 
$$ H_{A(t),\psi(t)}=\left\{ \begin{array}{c}
L_{A(t),\psi(t)}(\alpha,\phi) + G_{A(t),\psi(t)}(f) \\
G^*_{A(t),\psi(t)}(\alpha,\phi), \end{array} \right. $$
with
$$ L_{A(t),\psi(t)}(\alpha,\phi)=(*_3 d\alpha -2iIm(\psi(t),\phi),
\dirac_{A(t)} \phi + \alpha\cdot \psi(t)), $$
$$ G_{A(t),\psi(t)}(f)=(-df, e^{if}\psi(t)), $$
and $G_{A(t),\psi(t)}^*$ is the $L^2_{\delta}$ adjoint of
$G_{A(t),\psi(t)}$. Let $Q_{a_\infty}$ be the asymptotic operator of
$H_{A(t),\psi(t)}$, for each $t$, given by
$$ Q_{a_\infty}(a,\alpha,\beta)=(0,-i\bar\partial_{a_\infty}^* \beta,
i\bar\partial_{a_\infty}\alpha), $$
as in \cite{CMW}. 

Each moduli space 
$$ \M_{V\times \R}( [A,\psi],[\tilde A,\tilde
\psi], a_\infty), $$
for a fixed choice of $[A,\psi]$ and $[\tilde A,\tilde
\psi]$ in $\partial_\infty^{-1}(a_\infty)$ in $\M_V^*$, is a smooth
finite dimensional oriented manifold, of dimension given by the
spectral flow 
$$ SF(H_{A(t),\psi(t)})= Index( \partial_t + H_{A(t),\psi(t)}). $$
Thus, for the balanced energy case, we have virtual dimension
$$ \hbox{virtdim} \hat\M_{V\times \R}( [A,\psi],[\tilde A,\tilde
\psi], a_\infty) = SF(H_{A(t),\psi(t)}) -1. $$
An orientation of $\hat \M_{V\times \R}( [A,\psi],[\tilde A,\tilde
\psi], a_\infty)$ is obtained by considering the determinant line
bundle of the operator ${\cal L}_{\A,\Psi}=\partial_t + H_{A(t),\psi(t)}$.
\end{Pro}

\subsection{Admissible elements in $\hat \M_{V\times \R}(a_\infty)$}

Now consider $a(t)$ a regular parameterization of an arc in
$\partial_\infty(\M^*_V)$ inside  
$\chi_0(T^2,V)$, which satisfies
$$ a(-1)=a^- \ \ a(0)=a_\infty \ \ a(1)=a^+, $$
with given 
$$ a^- \in \partial_\infty(\M^*_V)\cap \{ u-v=1 \} $$
$$ a^+ \in \partial_\infty(\M^*_V)\cap \{ v=f'(u) \}. $$

Let us assume that the path $a(t)$ also satisfies the condition
$a(t)\neq \vartheta$, for all $t\in [-1,1]$, and that the path $a(t)$
avoids all the boundary points of $\partial_\infty(\M_V^*)$ on
$\chi(V)\backslash \vartheta$, cf. \cite{CMW}. 
By the analysis of \cite{CMW}, we know then that the fiber
$$ \partial_\infty^{-1}(a(t)) \subset \M^*_V $$
is a finite set of points, for each fixed $t$. Moreover, under the
current hypotheses, the set
$$ \cup_{t\in [-1,1]} \partial_\infty^{-1}(a(t)) $$
describes a cobordism between 
$$\partial_\infty^{-1}(a(-1)) \ \  \hbox{ and } \ \
\partial_\infty^{-1}(a(0))$$ and 
between $$ \partial_\infty^{-1}(a(0)) \ \  \hbox{ and } \ \
\partial_\infty^{-1}(a(1)). $$ 

\begin{figure}[ht]
\epsfig{file=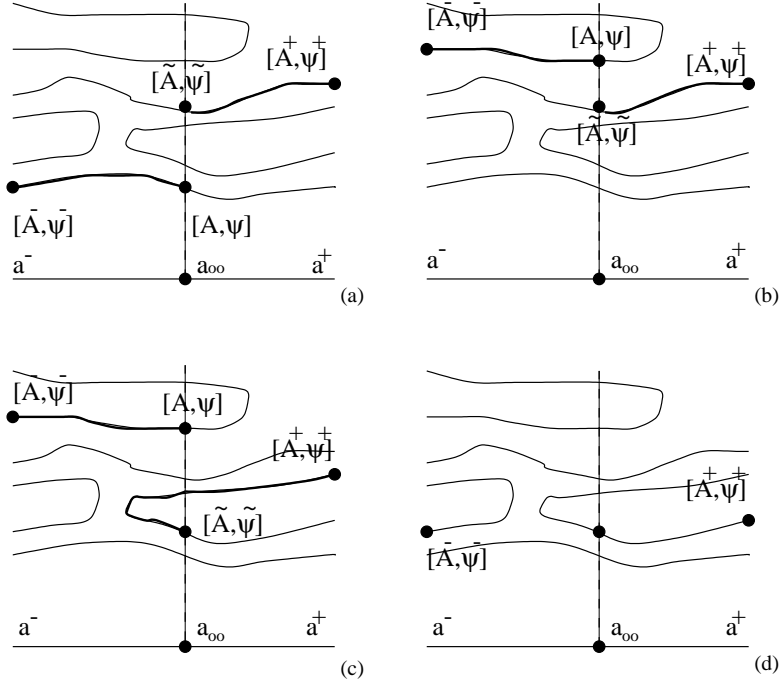,angle=270}
\caption{The paths $[A(t),\psi(t)]$ in $\M^*_V$ \label{figIIIpaths} }
\end{figure}

Now suppose given two assigned elements 
$$ [A^-,\psi^-] \in \partial_\infty^{-1}(a^-) $$
$$ [A^+,\psi^+] \in \partial_\infty^{-1}(a^+), $$
and elements
$$ [A,\psi], [\tilde A,\tilde \psi] \in
\partial_\infty^{-1}(a_\infty). $$
In the cobordism 
$$ \cup_{t\in [-1,1]} \partial_\infty^{-1}(a(t)), $$
there is at most one path $[A(t),\psi(t)]$ in $\M^*_V$, for $t\in
[-1,0]$ that satisfies $\partial_\infty [A(t),\psi(t)]=a(t)$ and
$[A(-1),\psi(-1)]=[A^-,\psi^-]$. Similarly, there is at most one
path $[A(t),\psi(t)]$ in $\M^*_V$, for $t\in [0,1]$ that satsifies
$\partial_\infty [A(t),\psi(t)]=a(t)$ and
$[A(1),\psi(1)]=[A^+,\psi^+]$, cf. 
the Figure \ref{figIIIpaths}. The case (d) of Figure
\ref{figIIIpaths} illustrates an example of a parameterization $a(t)$ and
a choice of $[A^-,\psi^-]$ and $[A^+,\psi^+]$ for which a path with
the desired properties does not exist. Such case does not arise as a
geometric limit. 

Thus, we can introduce the following notation, to distinguish which
choices of elements $([A,\psi],[\tilde A,\tilde \psi], a_\infty)$ can
arise as part of the geometric limits of solutions in
$\M^{W_1(r)}_\ell(a_1,a)$  
(or $\M^{W_0(r)}_k(a_1,a_0)$, or $\M_{Y(r)\times \R}(a,b)$).

\begin{Def}
We say that a triple $([A,\psi],[\tilde A,\tilde \psi], a_\infty)$ is
{\em admissible}, with respect to the endpoints $(a_1,a)$,
if the following conditions hold. The element
$a_\infty$ lies on a path component of $\partial_\infty(\M_V^*)$
connecting $a^-$ and $a^+$. Moreover, there exists a smooth regular
parameterization $a(t)$, for $t\in [-1,1]$ of the path in
$\partial_\infty(\M_V^*)$ connecting $a^-$ and $a^+$, such that
$a(0)=a_\infty$, and corresponding smooth paths $[A(t),\psi(t)]$ in
$\M_V^*$, for $t\in [-1,0)$ and $t\in (0,1]$, satisfying  
$\partial_\infty[A(t),\psi(t)]=a(t)$, and with
$$ [A(-1),\psi(-1)]=[A^-,\psi^-] \ \ \hbox{ and
} \ \ \lim_{t\to
0_-}[A(t),\psi(t)]=[A,\psi] $$
$$ \lim_{t\to 0_+}[A(t),\psi(t)]=[\tilde A,\tilde \psi] \ \ \hbox{ and
} \ \ [A(1),\psi(1)]=[A^+,\psi^+], $$
with $[A,\psi]$ and $[\tilde A,\tilde \psi]$ in
$\partial_\infty^{-1}(a_\infty)$ in $\M_V^*$ and
$$ a_1=[(A^-,\psi^-)\#(a^-,0)] $$
$$ a=[(A^+,\psi^+)\#(a^+,0)]. $$ 
\label{admissible}
\end{Def}

In the examples of Figure \ref{figIIIpaths}, the cases (a)-(b)
represent admissible elements, and case (c) is not an admissible
element, because the parameterization $a(t)=\partial_\infty
[A(t),\psi(t)]$ is not regular ($\partial_t a(t)= 0$ for some $t$),
and the case (d) is also non admissible because no path
$[A(t),\psi(t)]$ with the desired properties exists.
An element $([A,\psi],[\tilde A,\tilde \psi],
a_\infty)$ can appear as part of the geometric limits of solutions on
$W_1(r)$ (or $W_0(r)$, or $Y(r)\times \R$, etc.) only if it is
admissible with respect to the endpoints $(a_1,a)$ (or $(a,a_0)$, or
$(a,b)$).

Notice that, in general, if we consider different solutions in
$\M^{W_1(r)}_\ell(a_1,a)$ (or $\M^{W_0(r)}_k(a_1,a_0)$, or
$\M_{Y(r)\times \R}(a,b)$), these will give rise to geometric limits
with different parameterizations $a(t)$, and in general different
$a_\infty$. In the next subsection
we describe how to assemble the various geometric limits.

\subsection{Assembling the geometric limits}

Let us first return to the setting of \cite{MW2} and consider moduli
spaces of flowlines. We shall then generalize our statements to the
case of finite energy solutions on one of the surgery cobordisms.

Consider first the case of a zero-dimensional moduli space of flow lines,
$\hat\M_{Y\times\R,\mu}(a,b)$, with $a,b\in \M_{Y,\mu}^*$ satisfying
$\deg_{Y,\mu}(a)-\deg_{Y,\mu}(b)=1$. We recall the results on the
splitting of the spectral flow that we used in Part I \cite{CMW} in
order to compare the relative gradings, cf. \cite{DaKi}. 

We can describe the critical points $a$ and $b$ in $\M_{Y,\mu}^*$ as
$$ a=[(A^-,\psi^-)\#_r (a^-,0)] $$
$$ b=[(A^+,\psi^+)\#_r (a^+,0)]. $$
The relative grading is given by
$$ 1= \deg_{Y,\mu}(a)-\deg_{Y,\mu}(b) =
\frac{1}{r^2}SF_{Y(r)}(H_{(A_r(t),\psi_r(t))}), $$
with the notation as in \cite{CMW}, \cite{CLM}.
Then the spectral flow can be written
as 
\ba \begin{array}{l}\frac{1}{r^2}SF_{Y(r)}(H_{(A_r(t),\psi_r(t))})= \epsilon
SF_{V(r)}(H_{(A(t),\psi(t)),\tilde\ell_V(t)})+ \\[2mm]
\epsilon SF_{\nu(K)(r)}(H_{(a''(t),0),\tilde\ell_\nu(t)}) +
Maslov(\tilde\ell_V(t),\tilde\ell_\nu(t)) \end{array} 
\label{maslov:split} \na
Here we follow the same convention of \cite{CLM} regarding the definition
of the $\epsilon$-spectral flow.
The  boundary conditions are prescribed by assigning a choice
of Lagrangian subspaces $(\tilde\ell_V(t),\tilde\ell_\nu(t))$
in $H^1(T^2,\R)$. 

We consider the following
Lagrangians, which we already introduced in Part I, \cite{CMW}.
Let $\ell_\mu^*$ be the piecewise smooth Lagrangian submanifold of
$\chi_0(T^2,Y)$ described in Part I, \cite{CMW}, and let $\tilde
\ell_\mu$ be the path of Lagrangian subspaces of $H^1(T^2,\R)$ defined
as in Part I, 
\cite{CMW}, which is given by the tangent spaces of $\ell_\mu^*$ where
the latter is smooth, completed with a specific choice of paths at the
singular points, as discussed in \cite{CMW}. Let $\ell$ be the union
of the arcs in the Lagrangian submanifold with boundary
$\partial_\infty(\M_V^*)$ in $\chi_0(T^2,Y)$ that connect the points
$a^-$ and $a^+$. By $\partial_\infty(\M_V^*)\subset \chi_0(T^2,Y)$ we
mean the pullback of $\partial_\infty(\M_V^*) \subset \chi_0(T^2,V)$
under the covering map $\chi_0(T^2,Y)\to \chi_0(T^2,V)$.
Under the assumption that these arcs avoid the boundary points 
$\partial_\infty(\partial \M_V^*)$, we can consider regular
parameterizations $a'(t)$ and $a''(t)$ of the arcs of $\ell_\mu^*$ and
$\ell$, respectively,  connecting $a^-$ to $a^+$. We consider the
corresponding paths of Lagrangians $\tilde\ell(t)=T_{a'(t)}\ell$ and
$\tilde\ell_\mu(t)$, for $t\in [0,1]$. We can also assume that 
$\ell_\mu^*$ and $\ell$ intersect transversely.

Now recall that, in Part II, we identified the geometric limits of
flowlines in $\hat\M_{Y(r)\times\R}(a,b)$ with finite energy solutions
in $\hat\M_{V\times\R}([A,\psi],[\tilde A,\tilde \psi],a_\infty)$,
together with paths $[A(t),\psi(t)]\in \M_V^*$, and a holomorphic disk
$$ \Delta: (D^2,\partial D^2) \to (\chi_0(T^2,Y),\ell\cup \ell_\mu^*),
$$ 
which determines the regular parameterization $a(t)$ of the arcs of
$\ell$ and $\ell_\mu^*$, hence the admissible data
$([A,\psi],[\tilde A,\tilde \psi], a_\infty)$. Thus, we are interested in
understanding the space of inequivalent holomorphic disks $\Delta$,
which can appear in the geometric limits. 
We have the following result, see for instance
\cite{Fuk} pg.3. The set of equivalence classes under $Aut(D^2)\simeq
PSL(2,\R)$ of holomorphic disks
$$ \Delta: (D^2,\partial D^2) \to (\chi_0(T^2,Y),\ell\cup \ell_\mu^*),
$$ 
in a given homotopy class $\beta=[\Delta]$ in $\pi_2(\chi_0(T^2,Y), \ell\cup
\ell_\mu^*)$ has virtual dimension
$$ \mu(\beta)-1, $$
with $\mu(\beta)$ the Maslov class, defined as in
\cite{Fuk} pg. 3. The original formula given in \cite{Fuk} for this
virtual dimension is $ n +
\mu(\beta)-2$, where $2n$ is the real dimension of the ambient
symplectic manifold. In our case, we have $n=1$.
We can describe the Maslov class $\mu(\beta)$ in terms of the Maslov
index $Maslov(\tilde\ell(t),\tilde\ell_\mu(t))$ as follows.

\begin{Lem}
\label{maslov-maslov}
Consider the piecewise smooth Lagrangian submanifolds $\ell$ and
$\ell_\mu^*$ in $\chi_0(T^2,Y)$, defined as above. 
Consider all possible holomorphic disks
$$ \Delta: (D^2,\partial D^2) \to (\chi_0(T^2,Y),\ell\cup \ell_\mu^*),
$$ 
up to automorphisms $Aut(D^2)\simeq PSL(2,\R)$, which map the boundary
to arcs of the Lagrangians connecting two points $a^-$ and $a^+$ in
$\ell\cap \ell_\mu^*$, in the homotopy class $\beta$ specified by the
regular parameterizations $a'(t)$ and $a''(t)$ of 
these arcs of Lagrangians. Then the Maslov class satisfies 
$$ \mu(\beta) = Maslov(\tilde\ell(t),\tilde\ell_\mu(t))+1. $$
\end{Lem}

\proof
Here we follow the notation of Section 13 of \cite{CLM-Maslov}. First
let us define 
$$ \nu_{\pm}^{a^\pm}(T_{a^\pm}\ell_1,T_{a^\pm}\ell_2) $$
as the path of Lagrangian subspaces that connects the two Lagrangian
specified subspaces $(T_{a^+}\ell_1,T_{a^+}\ell_2)$, rotating 
in the positive or negative direction, according to sign. Recall that,
given a pair of Lagrangian paths
$f(t)=(\tilde\ell_1(t),\tilde\ell_2(t))$, we use the notation $\hat
f(t)=f(1-t)$, and 
$$ f_{flip}(t)=(\tilde\ell_2(t),\tilde\ell_1(t)). $$
Thus, given $f(t)=(\tilde\ell_1(t),\tilde\ell_2(t))$ with 
$\ell_1(0)=\ell_2(0)=a^-$ and $\ell_1(1)=\ell_2(1)=a^+$, we obtain
four possible loops $f_{\pm\pm}(t)$, given as the concatenation of
paths
$$ f_{ij}(t)= \tilde\ell_1(t) * \nu_{i}^{a^+}(t) *
\widehat{\tilde\ell_2(t)} * (\nu_{j}^{a^-}(t))_{flip}, $$
with $i,j=\pm$.

The definition of the Maslov class given in \cite{Fuk} pg. 3, for one
smooth Lagrangian submanifold $\ell$ inside a symplectic manifold $X$
and disks $\Delta: (D^2,\partial D^2) \to (X, \ell)$, coincides with
the notion of Maslov index for closed loops {\em \`a la} Floer
\cite{Floer}. Using the results of Section 13 of \cite{CLM-Maslov}, we
can relate this to the Maslov index, as described in the Lemma. Since
the arguments of Part II, \cite{MW2} actually ensure the existence of
a holomorphic disk on a slightly larger domain, which maps a subdomain
homeomorphic to $(D^2,\partial D^2)$ to a disk filling a region
bounded by arcs of the Lagrangians $\ell\cup \ell_\mu^*$ with
endpoints $a^-$ and $a^+$, we obtain that the Maslov class is computed
by 
$$ \mu(\beta)= \mu(f_{++}(t)), $$
with the notation as above, and $\mu(f_{++}(t))$ the Maslov index of a
loop of Lagrangians {\em \`a la} Floer, cf. Section 13 of
\cite{CLM-Maslov}, where we have 
$$ f(t)=(\tilde\ell(t),\tilde\ell_\mu(t)). $$ 
Now the result of Section 13 of \cite{CLM-Maslov} gives 
$$ \mu(f_{++}(t))= Maslov(\tilde\ell(t),\tilde\ell_\mu(t)) +1. $$

\endproof

Notice that we can think of these holomorphic disks equivalently as
equivalence classes of 
$$ \Delta: (D^2,\partial D^2)\to (\chi_0(T^2,Y),\ell\cup \ell_\mu^*)
$$
that map $\Delta(-1)=a^-$ and $\Delta(+1)=a^+$, modulo the subgroup of
$Aut(D^2)$ that fixes the points $\pm 1$, or then again, equivalently,
as holomorphic maps $\Delta$ of an infinite strip
$$  \R \times [0,1] \to \chi_0(T^2,Y) $$
that map the boundary $\R \times \{ 0,1 \}$ to arcs of the
Lagrangians $\ell$ and $\ell_\mu^*$, and with asymptotic values along
$t\in \R$, as $t\to \mp \infty$, equal to $a^-$ and $a^+$. In this
case, we consider classes modulo the action of $\R$ by
reparameterizations.

Thus, the solutions in the zero-dimensional moduli space
$\hat\M_{Y\times\R,\mu}(a,b)$ are obtained by gluing solutions in a
zero-dimensional moduli space $$\hat \M_{V\times
\R}([A,\psi],[\tilde A,\tilde \psi], a_\infty)$$ with a holomorphic disk
$$ \Delta: (D^2,\partial D^2) \to (\chi_0(T^2,Y),\ell\cup \ell_\mu^*), $$
in a fixed homotopy class $\beta=[\Delta]$ in $\pi_2(\chi_0(T^2,Y),
\ell\cup \ell_\mu^*)$. The admissible data $([A,\psi],[\tilde A,\tilde \psi],
a_\infty)$ is determined by the holomorphic disk $\Delta$, as we
discuss in the following Lemma.

\begin{Lem}
Suppose given a zero-dimensional moduli space
$\hat\M_{Y\times\R,\mu}(a,b)$. 
The space of pre-glued solutions, obtained by pasting together the
geometric limits of solutions in $\hat\M_{Y\times\R,\mu}(a,b)$, as
described in \cite{MW2}, is given by
$$ \bigcup_{\Delta\in \Xi} \hat \M_{V\times \R}([A,\psi],[\tilde A,\tilde
\psi], a_\infty). $$
where $\Delta$ is a choice of one particular representative in each
equivalence class in the set $\Xi$ of classes of holomorphic disks
$$ \Delta: (D^2,\partial D^2) \to (\chi_0(T^2,Y),\ell\cup \ell_\mu^*),
$$
in the fixed homotopy class $[\Delta]$ in $\pi_2(\chi_0(T^2,Y), \ell\cup
\ell_\mu^*)$, modulo the action of $Aut(D^2)\simeq PSL(2,\R)$. 
Each such representative $\Delta$ determines the corresponding
admissible element 
$([A,\psi],[\tilde A,\tilde \psi], a_\infty)$.
\label{pre-glue:V}
\end{Lem}

\proof
The statement of the Lemma follows from the previous discussion. In
order to see the dependence of the data $([A,\psi],[\tilde A,\tilde
\psi], a_\infty)$ on the holomorphic disk $\Delta$, recall that the
choice of $\Delta$ in particular fixes the parameterization $a(t)$ of
the arc of $\ell$ connecting $a^-$ and $a^+$. This determines the point
$a_\infty=a(0)$. However, the choice of the parameterization $a(t)$ also
determines the choice of the points $[A,\psi]$ and $[\tilde A,\tilde
\psi]$ in $\partial_\infty^{-1}(a_\infty)$.
Notice that each solution $(\A(r),\Psi(r))$, representing an element
in $\hat\M_{Y(r)\times\R,\mu}(a,b)$, determines a holomorphic disk
$$ \Delta=\Delta(\A(r),\Psi(r)) $$
in the geometric limits. Each such disk determines a possibly
different parameterization $a(t)$, hence different admissible data 
$$ ([A,\psi],[\tilde A,\tilde
\psi], a_\infty). $$
Thus, the set of pre-glued solutions can be written as
$$ \bigcup_{\Delta} \hat \M_{V\times \R}([A,\psi],[\tilde A,\tilde
\psi], a_\infty), $$
where $\Delta$ varies in the set of holomorphic disks
$$ \Delta: (D^2,\partial D^2) \to (\chi_0(T^2,Y),\ell\cup \ell_\mu^*),
$$
in the fixed homotopy class $[\Delta]$ in $\pi_2(\chi_0(T^2,Y), \ell\cup
\ell_\mu^*)$. However, we only count the disks $\Delta$ up to
$Aut(D^2)\simeq PSL(2,\R)$. 

For geometric limits of a zero-dimensional moduli space,
the set of such holomorphic disks, up to $Aut(D^2)$ is also a
zero-dimensional moduli space. This follows by the result of Lemma
\ref{maslov-maslov}, in the case where we have
$$Maslov(\tilde\ell(t),\tilde\ell_\mu(t))=0.$$ In fact, in our case, with
disks in $\chi_0(T^2,Y)$, the homotopy class is fixed by the choice of
the arcs of the Lagrangians $\ell$ and $\ell_\mu^*$ connecting $a^-$
and $a^+$, and the representatives in each class only differ by
reparameterizations of the boundary. A choice of a representative
$\Delta$ in each class identifies uniquely a corresponding moduli
space 
$$ \hat \M_{V\times \R}([A,\psi],[\tilde A,\tilde
\psi], a_\infty). $$

Thus, with the notation as above, we obtain
$$ n_{Y,\mu}(a,b)=\sum_{\Delta} \# \hat \M_{V\times
\R}([A,\psi],[\tilde A,\tilde \psi], a_\infty). $$ 

\endproof

Thus, we can rephrase the gluing theorem of \cite{MW2} as the
following statement. Let us introduce the following notation:
$$ \hat\M_{V\times \R}([A^-,\psi^-],[A^+,\psi^+]):=
\bigcup_{\Delta} \hat \M_{V\times \R}([A,\psi],[\tilde A,\tilde
\psi], a_\infty), $$
with the admissible data $([A,\psi],[\tilde A,\tilde
\psi], a_\infty)$ and the corresponding holomorphic disk
$\Delta$ as in Lemma \ref{pre-glue:V}.

\begin{Pro}
The gluing map gives an orientation preserving diffeomorphism
$$ \#: \hat\M_{V\times \R}([A^-,\psi^-],[A^+,\psi^+]) \cong
\hat\M_{Y\times \R}(a,b), $$
with 
$$ a=[(A^-,\psi^-)\#_{a^-}(a^-,0)] $$
$$ b=[(A^+,\psi^+)\#_{a^+}(a^+,0)]. $$
\label{glue:V}
\end{Pro}

Now we can rephrase the result in the case of moduli spaces of finite
energy solutions of the Seiberg--Witten equations on one of the
surgery cobordisms. In this case, we can assemble the geometric limits
in a similar way, to obtain spaces
\ba \M_{V\times \R}([A^-,\psi^-],[A^+(\epsilon),\psi^+(\epsilon)]):=
 \bigcup_{\Delta_1^\epsilon} \M_{V\times
\R}([A,\psi]_1^\epsilon,[\tilde A,\tilde \psi]_1^\epsilon, 
a_\infty^1(\epsilon)), \label{def:M:V:W1} \na 
in the case of $W_1$, with the admissible data
$$([A,\psi]_1^\epsilon,[\tilde A,\tilde \psi]_1^\epsilon,
a_\infty^1(\epsilon))$$ 
determined by the holomorphic triangle $\Delta_1^\epsilon$,
or
\ba \M_{V\times \R}([A^-(\epsilon),\psi^-(\epsilon)],[A^+,\psi^+]):=
\bigcup_{\Delta_0^\epsilon}  \M_{V\times
\R}([A,\psi]_0^\epsilon,[\tilde A,\tilde \psi]_0^\epsilon, 
a_\infty^0(\epsilon)), \label{def:M:V:W0} \na in the case of $W_0$,
with the admissible data
$$ ([A,\psi]_0^\epsilon,[\tilde A,\tilde \psi]_0^\epsilon,
a_\infty^0(\epsilon)) $$
determined by the holomorphic triangle $\Delta^\epsilon_0$. 
Here the $\Delta_i^\epsilon$ vary in the set of inequivalent
holomorphic triangles in 
$H^1(T^2,\R)$, namely 
$$ \Delta_1^\epsilon: (D^2,\partial D^2) \to (H^1(T^2,\R),\ell \cup
\ell_1 \cup \ell_\mu), $$
in the case of $W_1$, or
$$ \Delta_0^\epsilon: (D^2,\partial D^2) \to (H^1(T^2,\R),\ell \cup
\ell_\mu \cup_k \ell_k), $$ in the case of $W_0$,
with $\ell_1=\{ u-v=1 \}$, $\ell_\mu=\{ v=f'(u) \}$, and $\ell_k=\{ u=
2k \}$ or $\{ u=\eta \}$ in the case of $\ell_0$, as in \cite{CMW}. 
The image of each $\Delta_i^\epsilon$ describes a triangle in
$H^1(T^2,\R)$ with 
vertices $\{ a^-, \vartheta_i, a^+(\epsilon) \}$ or
$\{ a^-(\epsilon), \vartheta_i, a^+ \}$ and sides along the Lagrangians,
as specified. Here the points $\vartheta_i$ are the intersection of
the lines $\ell_1$ and $\ell_\mu$ for $W_1$ and $\ell_k$ and
$\ell_\mu$ for $W_0$, that is, the restriction to $T^2=\partial
\nu=\partial \tilde \nu$ of the unique reducible point $\theta_{S^3}$
at the puncture in the cobordism.

Thus, we obtain the following result on the gluing theorem for the
moduli spaces $\M_\ell^{W_1}(a_1,a)$, or for the minimal energy
component of $\M_{k}^{W_0}(a,a_0)$. We state the result in the case of
$W_1$.

\begin{Lem}
Suppose given a pair $a_1 \in \M_{Y_1}$, and $a \in \M_{Y,\mu}$.
Suppose that we have the decomposition
$$ a_1=[(A^-,\psi^-)\#_{a^-}(a^-,0)] $$
$$ a=[(A^+(\epsilon),\psi^+(\epsilon))
\#_{a^+(\epsilon)}(a^+(\epsilon),0)].  $$ 
Then the gluing map gives an orientation preserving diffeomorphism 
$$ \#: \M_{V\times
\R}([A^-,\psi^-],[A^+(\epsilon),\psi^+(\epsilon)]) \to
\M_\ell^{W_1}(a_1,a), $$
where the first moduli space is defined as in (\ref{def:M:V:W1}), with
the union over inequivalent holomorphic triangles $\Delta$ with vertices
$\{ a^-,\vartheta_1, a^+(\epsilon) \}$ and sides along the union of
Lagrangians $\ell \cup \ell_1 \cup \ell_\mu$, with $\ell$ defined by
the asymptotic values $\partial_\infty(\M_V^*)$. Different choices of the
$\spinc$-structure $\s_\ell$ correspond to moduli spaces of different
dimension. If $\M^{W_1}_\ell(a_1,a)$ is non-empty, then, under the
gluing map we obtain 
$$ \iota(W_1, \s_\ell, a_1,a)= \mu([\Delta])-1 + virtdim \M_{V\times
\R}( [A,\psi],[\tilde A,\tilde \psi], a_\infty ). $$
\label{diff:V:W_1}
\end{Lem}

Again we can observe that each holomorphic triangle fixes a
parameterization $a(t)$, 
hence a choice of the admissible data $([A,\psi],[\tilde A,\tilde
\psi], a_\infty)$. The argument then proceeds as in the case of flow
lines. In the case of the moduli spaces $\M_{k}^{W_0}(a,a_0)$, the minimal
energy condition ensures that the path $a(t)$ is contractible in
$\chi_0(T^2,V)$, so that the path along the union of Lagrangians can
be filled by a holomorphic triangle of minimal energy.

We now return to our analysis of the sequence
$$ 0\to C_*(Y_1)\stackrel{w^1_*}{\to} C_*(Y,\mu)\stackrel{w^0_*}{\to}
\bigoplus_{k} C_{(*)} (Y_0,\s_{k})\to 0. $$
In the next subsection, we use the results obtained in this section on
the moduli spaces 
$\M_{V\times \R}( [A^-,\psi^-],[A^+,\psi^+])$ to prove exactness
in the middle term.

\subsection{The relation $w_*^0\circ w_*^1=0$}

By Lemma \ref{comp:rel}, we can show that the composite map
$w_*^0\circ w_*^1$ is trivial by proving the following result. 

\begin{The}\label{W_1:W_0} 
For small enough $\epsilon$ and large $r\ge r_0$, there 
is an orientation reversing  diffeomorphism  
\[   
\M_\ell ^{W_1(r)}\bigl(a_i^{(1)}, a_j^{(0)}(\epsilon) \bigr) 
\cong \M_k^{W_0(r),\  (0)} \bigl(a_i^{(1)}(\epsilon), a_j^{(0)}\bigr), 
\] 
where $\ell$ is the unique positive $\spinc$ structure on $W_1(r)$ such that 
the moduli space $\M_\ell ^{W_1(r)}\bigl(a_i^{(1)},
a_j^{(0)}(\epsilon) \bigr)$  
is zero dimensional, and  
$\M_k^{W_0(r),\  (0)} \bigl(a_i^{(1)}(\epsilon), a_j^{(0)}\bigr)  $ 
is the zero-dimensional components of 
 $\M_k^{W_0(r) } \bigl(a_i^{(1)}(\epsilon), a_j^{(0)}\bigr)$. 
\end{The} 

\proof 
By our assumptions, the moduli spaces $\M_\ell ^{W_1(r)}\bigl(a_i^{(1)},
a_j^{(0)}(\epsilon) \bigr)$  
and $\M_k^{W_0(r),\  (0)} \bigl(a_i^{(1)}(\epsilon), a_j^{(0)}\bigr)$ are 
smooth, compact, oriented 0-dimensional manifolds.  
We write the critical points $a_i^{(1)}$, 
$a_j^{(0)}(\epsilon)$ and $a_i^{(1)}(\epsilon), a_j^{(0)}$ according
to the following decomposition:  
\[ 
a_i^{(1)} = [(A_i(r), \psi_i(r))] = [(A_i^- , \psi_i^-) 
\#_{a_{\infty, i}^-} (a_{\infty, i}^-, 0)] 
\] 
on $Y_1(r)$,  
\[\begin{array}{l} 
a_i^{(1)}(\epsilon) = j ([(A_i(r), \psi_i(r))])  
= [(A_i^-(\epsilon), \psi_i^-(\epsilon)) \#_{a_{\infty, i}^-(\epsilon)}
(a_{\infty, i}^-(\epsilon), 0)]\\[2mm]  
a_j^{(0)}(\epsilon) = [(A_j^+(\epsilon) , \psi_j^+(\epsilon)) 
\#_{a_{\infty, j}^+(\epsilon)} (a_{\infty, j}^+(\epsilon), 0)] 
\end{array} 
\] 
on $Y(r)$, and  
\[ 
a_j^{(0)} =  [( A_j^+ ,\psi_j^+) \#^r_{ a_{\infty, j}^+} (
a_{\infty, j}^+, 0)]  
\] 
on $Y_0(r)$. In the limit $\epsilon \to 0$, the asymptotic values 
satisfy
$$ \lim_{\epsilon \to 0} a_{\infty, i}^-(\epsilon) =a_{\infty, i}^- $$
$$ \lim_{\epsilon \to 0} a_{\infty, j}^+(\epsilon) =a_{\infty, j}^+,
$$
in $H^1(T^2,i\R)$ and
$$ \lim_{\epsilon \to 0} [A_i^-(\epsilon), \psi_i^-(\epsilon)]=[A_i^- ,
\psi_i^-] $$
$$ \lim_{\epsilon \to 0} [A_j^+(\epsilon) , \psi_j^+(\epsilon)] =
[A_j^+ ,\psi_j^+] $$ 
in $\M_V^*$.

We now apply  
the geometric limits results of Proposition \ref{geom:lim} to study 
the moduli spaces $\M_\ell^{W_1(r)}\bigl(a_i^{(1)},
a_j^{(0)}(\epsilon)\bigr)$  
and $\M_k^{W_0(r),\  (0)} \bigl(a_i^{(1)}(\epsilon), a_j^{(0)}\bigr)$. 
 
Let $[\A_1(r), \Psi_1(r)]$ be solutions in  
$\M_\ell^{W_1(r)}\bigl(a_i^{(1)}, a_j^{(0)}(\epsilon)\bigr)$, for 
sufficiently large $r \geq r_0$. Then these solutions define geometric 
limits as recalled in Remark \ref{geom:lim:V}, cf. Proposition
\ref{geom:lim}.

For all sufficiently small $\epsilon$, the geometric
limits listed in Proposition \ref{geom:lim} can be grouped together to
give a holomorphic triangle 
$$ \Delta_1^\epsilon: (D^2,\partial D^2)\to
(H^1(T^2,\R),\ell\cup\ell_1\cup\ell_\mu) $$   
with vertices $\{ a_{\infty,i}^- , \vartheta_1,
a_{\infty,j}^+(\epsilon) \}$ and sides along the Lagrangians
$\ell_1=\{ v-u=1 \}$, $\ell=\partial_{\infty} \M_V^*$, and $\ell_\mu=\{
v=f'(u) \}$, and an element in a moduli
space 
$$\M_{V\times\R}([A,\psi]_1^\epsilon,[\tilde A,\tilde
\psi]_1^\epsilon,a_\infty^1(\epsilon)), $$
for some admissible data 
$$([A,\psi]_1^\epsilon, [\tilde A,\tilde \psi]_1^\epsilon,
a_\infty^1(\epsilon))$$ 
determined by the holomorphic triangle $\Delta_1^\epsilon$. 

We use the notation  
$\Delta_1^\epsilon=\Delta^\epsilon(\A_1(r), \Psi_1(r))$, for any such 
holomorphic triangle determined by a family of solutions
$(\A_1(r),\Psi_1(r))$ representing an element in 
$$ \M_\ell^{W_1(r)}\bigl(a_i^{(1)}, a_j^{(0)}(\epsilon)\bigr). $$
The orientation of this triangle is the same as the orientation of the 
region of non-uniform convergence. 
The orientation of 
$$\hat\M_{V\times\R}([A,\psi]_1^\epsilon,[\tilde A,\tilde
\psi]_1^\epsilon,a_\infty^1(\epsilon))$$ 
at a solution $[\A',\Psi']^\epsilon$ is given by the
determinant line bundle of the linearization ${\cal
L}_{(\A',\Psi'^\epsilon)}$, 
as discussed previously. 

Now we can proceed to compare the geometric limits of solutions in
$\M_\ell^{W_1(r)}\bigl(a_i^{(1)}, a_j^{(0)}(\epsilon)\bigr)$ 
and $\M_k^{W_0(r),\  (0)} \bigl(a_i^{(1)}(\epsilon), a_j^{(0)}\bigr)$.

Recall that, as discussed previously (cf. also \cite{MW2}), we can
pre-glue the geometric limits to form an approximate monopole on 
$W_1(r)$, for sufficiently large $r\ge
r_0$. That is, as discussed in the previous subsection, a
holomorphic triangle $\Delta$ determines the regular parameterizations
of the arcs of Lagrangians, 
hence the admissible data $([A,\psi],[\tilde A,\tilde
\psi],a_\infty)$. Therefore, the space of pre-glued monopoles can be
identified with the moduli space
$$  \M_{V\times\R}([A_i^-,\psi_i^-],
[A_j^+(\epsilon),\psi_j^+(\epsilon)])= \bigcup_{\Delta_1^\epsilon} 
\M_{V\times\R}([A,\psi]_1^\epsilon,[\tilde A,\tilde
\psi]_1^\epsilon,a_\infty^1(\epsilon)),  $$ 
with the admissible data specified by 
$$ \Delta_1^\epsilon=\Delta^\epsilon (\A_1(r),\Psi_1(r)), $$
ranging over the set of inequivalent holomorphic triangles in
$H^1(T^2,\R)$, with 
vertices $\{ a_{\infty,i}^-, \vartheta_1, a_{\infty, j}^+(\epsilon) \}$
and sides along arcs of the Lagrangians $\ell$, $\ell_1$, and
$\ell_\mu$ connecting these points. The orientation on the moduli space
$\hat\M_{V\times\R}([A_i^-,\psi_i^-], 
[A_j^+(\epsilon),\psi_j^+(\epsilon)])$ is the product orientation
of the pairs $$([\A',\Psi']^\epsilon, \Delta_1^\epsilon)$$ with 
$[\A',\Psi']^\epsilon$ an element of the moduli space
$$\hat\M_{V\times\R}([A,\psi]_1^\epsilon,[\tilde A,\tilde
\psi]_1^\epsilon,a_\infty^1(\epsilon))$$
determined by the corresponding triangle $\Delta_1^\epsilon$.
The gluing map gives an orientation preserving diffeomorphism
\ba 
\M_\ell^{W_1(r)}\bigl(a_i^{(1)}, a_j^{(0)}(\epsilon)\bigr) 
\cong \M_{V\times \R} \bigl([A_i^-,\psi_i^-],[
A_j^+(\epsilon),\psi_j^+(\epsilon)]\bigr), 
\label{W_1:V} 
\na 
as in Lemma \ref{diff:V:W_1} in the previous
subsection.
Moreover, by an analogous argument, we have a similar result for the
zero-dimensional component 
$ \M_k^{W_0(r),\  (0)} \bigl(a_i^{(1)}(\epsilon), a_j^{(0)}\bigr)$.

\noindent{\bf Claim:} For $ \M_k^{W_0(r),\  (0)}
\bigl(a_i^{(1)}(\epsilon), a_j^{(0)}\bigr)$, the zero-dimensional
component of $\M_k^{W_0(r)} \bigl(a_i^{(1)}(\epsilon),
a_j^{(0)}\bigr)$, there is no energy loss in the process of 
stretching along the $r$-direction. In this case, we have a similar
orientation preserving diffeomorphism  
\ba 
\M_k^{W_0(r),\  (0)} \bigl(a_i^{(1)}(\epsilon), a_j^{(0)}\bigr) 
\cong \hat \M_{V(\infty)\times \R}
\bigl([A_i^-(\epsilon),\psi_i^-(\epsilon)],[A_j^+,\psi_j^+]
\bigr). 
\label{W_0:V} \na
Here the moduli space of pre-glued solutions is given by
$$  \M_{V\times \R}
\bigl([A_i^-(\epsilon),\psi_i^-(\epsilon)],[A_j^+,\psi_j^+]=
 \bigcup_{\Delta^\epsilon_0} \M_{V\times
\R}([A,\psi]_0^\epsilon,[\tilde 
A,\tilde \psi]_0^\epsilon, a_\infty^0(\epsilon)), $$
with the admissible data determined by 
$$ \Delta_0^\epsilon=  \Delta^\epsilon(\A_0(r),\Psi_0(r)), $$
ranging over the set of inequivalent holomorphic triangles in
$H^1(T^2,\R)$, with vertices
$\{ a_{\infty,i}^-(\epsilon),\vartheta_0, a_{\infty,j}^+ \}$ and sides
along arcs of the Lagrangians $\ell$, $\ell_\mu$, and $\cup_k \ell_k$
connecting these points.

\noindent{\bf Claim:} As $\epsilon \to 0$, we can identify the moduli
spaces 
$$ \M_{V\times \R} \bigl([A_i^-,\psi_i^-],[
A_j^+(\epsilon),\psi_j^+(\epsilon)]\bigr) $$
and
$$ \M_{V\times \R}
\bigl([A_i^-(\epsilon),\psi_i^-(\epsilon)],[ A_j^+,\psi_j^+]
\bigr) $$
as sets of points.

In fact, if the triangles $\Delta^\epsilon_1$ and $\Delta^\epsilon_0$
are non-degenerate, for $\epsilon\to 0$, then these 
holomorphic triangles reduce to  
holomorphic triangles
$$ \Delta: (D^2,\partial D^2) \to (H^1(T^2,\R),
\ell\cup\ell_1\cup_k\ell_k), $$
in the same homotopy class in $\pi_2(H^1(T^2,\R),
\ell\cup\ell_1\cup_k\ell_k)$, with 
sides along arcs of the Lagrangians $\ell\cup\ell_1\cup_k\ell_k$
connecting the points $a_{\infty,i}^-$, $a_{\infty,j}^+$ and the
intersection point between the lines $\ell_1$ and $\ell_k$.
Notice that the set theoretic difference between the triangles 
$\Delta^\epsilon_1$ and $\Delta^\epsilon_0$ shrinks to zero size as
$\epsilon \to 0$, hence they define the same triangle $\Delta$ as a
limit. Again we can count
these up to automorphisms of $D^2$.
Thus, in the limit as $\epsilon\to 0$ we can identify both moduli
spaces with 
$$ \M_{V\times \R}
\bigl([A_i^-,\psi_i^-],[ A_j^+,\psi_j^+])=  \bigcup_{\Delta}
\M_{V\times\R}([A,\psi],[\tilde A,\tilde \psi], a_\infty), $$
with the admissible data determined by $\Delta$.
Notice that the parameterizations $a(t)$ determined by the holomorphic
triangles $\Delta$ are sufficiently close in the ${\cal C}^\infty$
topology to the parameterizations $a^1(t)$ and $a^0(t)$ determined by
the $\Delta_1^\epsilon$ 
and $\Delta_0^\epsilon$, for sufficiently small $\epsilon$, and the
corresponding admissible data $([A,\psi],[\tilde A,\tilde \psi],
a_\infty)$ are also sufficiently close to the admissible data
$$ ([A,\psi]_1^\epsilon,[\tilde A,\tilde \psi]_1^\epsilon,
a_\infty^1(\epsilon)) \ \ \hbox{ and } \ \
([A,\psi]_0^\epsilon,[\tilde A,\tilde \psi]_0^\epsilon,
a_\infty^0(\epsilon)). $$
Thus, for sufficiently close data,
the arguments of the previous subsections show that the moduli spaces 
$$  \begin{array}{l} \M_{V\times\R}([A,\psi],[\tilde
A,\tilde \psi], a_\infty) \cong  \\[2mm]
$$ \M_{V\times\R}([A,\psi]_1^\epsilon,[\tilde A,\tilde
\psi]_1^\epsilon, a_\infty^1(\epsilon)) \cong \\[2mm]
$$  \M_{V\times\R} ([A,\psi]_0^\epsilon,[\tilde A,\tilde
\psi]_0^\epsilon, a_\infty^0(\epsilon)) \end{array} $$
can also be identified.  

Thus, we have obtained that the moduli spaces
$\M_\ell^{W_1(r)}\bigl(a_i^{(1)}, a_j^{(0)}(\epsilon)\bigr)$ and
$\M_k^{W_0(r),\  (0)} \bigl(a_i^{(1)}(\epsilon), a_j^{(0)}\bigr)$ can
be identified as sets of points. This implies that the composite map
satisfies $w_*^0\circ w_*^1 \equiv 0 \, (mod \, 2)$. In order to
obtain the result over the integers, we need to compare the
orientations of these moduli spaces as they appear in the decomposition of
$\M_{\ell,k}^{W(r)}\bigl(a_i^{(1)},a_j^{(0)})$ in Lemma
\ref{split:W}.

\noindent{\bf Claim:} In the decomposition
$$ \M_{\ell,k}^{W(r)}\bigl(a_i^{(1)},a_j^{(0)}) \cong
\M_\ell^{W_1(r)}\bigl(a_i^{(1)}, a_j^{(0)}(\epsilon)\bigr) 
\cup
\M_k^{W_0(r),\  (0)} \bigl(a_i^{(1)}(\epsilon), a_j^{(0)}\bigr), $$
as in Lemma \ref{split:W}, the moduli spaces 
$$\M_\ell^{W_1(r)}\bigl(a_i^{(1)}, a_j^{(0)}(\epsilon)\bigr)\ \hbox{
and } \ \M_k^{W_0(r),\  (0)} \bigl(a_i^{(1)}(\epsilon),
a_j^{(0)}\bigr)$$ have  
opposite orientations.

In fact, it is sufficient to notice that, when we glue the punctured
cobordisms $\hat W_1(r)$ and $\hat W_0(r)$ along 
the common boundary $Y$, to obtain
\[  
 \hat W (r) = \hat W_1(r) \cup_{Y(r)} \hat W_0(r), 
\] 
there is a long product region from a solid torus in $Q_1$ to 
a solid torus in $Q_0$.
Thus, the orientations of the two rescaled  
regions of non-uniform convergence, that is, of the domains of the two
holomorphic triangles
$$\Delta^\epsilon_1=\Delta^\epsilon(\A_1(r),\Psi_1(r))\  \hbox{ 
and } \  \Delta^\epsilon_0=\Delta^\epsilon(\A_0(r),\Psi_0(r))$$ 
have opposite orientations. This  implies that, when regarded as
solutions in 
$$\M_{\ell,k}^{W(r)}\bigl(a_i^{(1)},a_j^{(0)}),$$ one of the triangles
is antiholomorphic (holomorphic up to a change of orientation). In
fact, the triangles in $H^1(T^2,\R)$ have the same orientation,
whereas the rescaled regions of non-uniform convergence for the
triangles on $\hat W_1$ and $\hat W_0$ differ by a change of
orientations when compared inside the composite cobordism $W$. 

In other words, we obtain an orientation reversing diffeomorphism 
\[ 
\M_\ell ^{W_1(r)}\bigl(a_i^{(1)}, a_j^{(0)}(\epsilon) \bigr) 
\cong \M_k^{W_0(r),\  (0)} \bigl(a_i^{(1)}(\epsilon), a_j^{(0)}\bigr), 
\] 
when both moduli spaces are identified with solutions in
$$\M_{\ell,k}^{W(r)}\bigl(a_i^{(1)},a_j^{(0)}). $$

\endproof 

We now derive from Theorem \ref{W_1:W_0} the third fundamental
relation. 

\begin{Cor} 
Let $a_i^{(1)}$ and $a_j^{(0)}$ be  critical points on  
$Y_1(r) $ and $Y_0(r)$ respectively with  
\[ 
deg_{Y_1}(a_i^{(1)}) = deg_{Y_0, \s_k}(a_j^{(0)}), 
\] 
and let $a_i^{(1)}(\epsilon)$ and $a_j^{(0)}(\epsilon)$ be the 
corresponding critical points on $Y(r)$, for sufficiently small $\epsilon$ 
and large $r\ge r_0$, then we have 
\[ 
   N_\ell^{W_1(r)}(a_i^{(1)}, a_j^{(0)}(\epsilon)) =
-N_k^{W_0(r)}(a_i^{(1)}(\epsilon),  
a_j^{(0)}).  
\] 
\label{N:-N}
\end{Cor}

\subsection{The inclusion $Ker(w_*^0)\subset Im(w_*^1)$}

Recall that, for a fixed pair of asymptotic values
$(a^{(1)}_i,a^{(0)}_j)$, there exists a unique non-negative $\spinc$
structure $\s_\ell$ satisfying
$$ \iota(\s_{\ell},W_1,a^{(1)}_i,a^{(0)}_j)=0. $$
We denote the corresponding $\ell$ with $\ell_{ij}$. The coefficient
of the map $w_*^1$ is then given by
$$N^{W_1}_{\ell_{ij}}(a^{(1)}_i,a^{(0)}_j).$$ 

Using the results of Section 6, we can write explicitly the condition
that a Floer chain on $Y$ is in the image of $w_*^1$.

\begin{Lem}
Suppose given a Floer chain $a\in C_*(Y,\mu)$,
$$ a= \sum_{i=1}^m n_i a^{(1)}_i + \sum_{j=m+1}^n n_j a^{(0)}_j, $$
with coefficients $n_i$ and $n_j$ in $\Z$.
If we assume that the parameter $\epsilon >0$ in the surgery
perturbation is sufficiently small, then the chain $a\in C_*(Y,\mu)$
is in the image $Im(w_*^1)$ under the morphism 
$$ w_*^1 : C_*(Y_1) \to C_*(Y,\mu) $$
if and only if the coefficients satisfy the relation
$$ n_j= \sum_{i=1}^m n_i  N^{W_1}_{\ell_{ij}}(a^{(1)}_i,a^{(0)}_j), $$
for all $j=m+1,\ldots, n$. 
\label{image}
\end{Lem}

\proof Assume that the element
$$ a= \sum_{i=1}^m n_i a^{(1)}_i + \sum_{j=m+1}^n n_j a^{(0)}_j $$
is in $Im(w_*^1)$. Then it satisfies
$$ a = w_*^1( \sum_{i=1}^m p_i  a^{(1)}_i), $$
for some coefficients $p_i \in \Z$.
For sufficiently small $\epsilon>0$, we know from the results of
Section 6 that the relation
$$ \la a^{(1)}_i, w_*^1 (a^{(1)}_k) \ra = \delta_{ik} $$
holds. This gives
$$ a= \sum_{i=1}^m p_i  a^{(1)}_i + \sum_{j=m+1}^n \sum_{i=1}^m
p_i N^{W_1}_{\ell_{ij}}(a^{(1)}_i,a^{(0)}_j). $$

Thus, we obtain
$$ n_i =p_i \ \hbox{ for all } \ i=1,\ldots, m $$
and
$$ n_j= \sum_{i=1}^m p_i N^{W_1}_{\ell_{ij}}(a^{(1)}_i,a^{(0)}_j). $$

\endproof

On the other hand, we can also write explicitly the condition that a
Floer chain on $Y$ is in the kernel of $w_*^0$, using the expression
for the map $w_*^0$ for sufficiently small $\epsilon >0$ in the
surgery perturbation. 

\begin{Lem}
Suppose given a Floer chain $a\in C_*(Y,\mu)$,
$$ a= \sum_{i=1}^m n_i a^{(1)}_i + \sum_{j=m+1}^n n_j a^{(0)}_j, $$
with coefficients $n_i$ and $n_j$ in $\Z$.
For sufficiently small $\epsilon >0$ in the
surgery perturbation, we have that $a\in C_*(Y,\mu)$ is in 
the kernel $Ker(w_*^0)$ if and only if the coefficients satisfy
$$ n_j = - \sum_{i=1}^m n_i N^{W_0}_{k_j}(a^{(1)}_i,a^{(0)}_j), $$  
for all $j=m+1,\ldots, n$, with
$$ a^{(0)}_j \in \M_{Y_0}(\s_{k_j}). $$
\label{kernel}
\end{Lem}

\proof  Assume that the element
$$ a= \sum_{i=1}^m n_i a^{(1)}_i + \sum_{j=m+1}^n n_j a^{(0)}_j $$
is in $Ker(w_*^0)$. Then it satisfies
$$ \left( \sum_{i=1}^m \sum_{j=m+1}^n n_i N^{W_0}_{k_j}(a^{(1)}_i,a^{(0)}_j)
 + \sum_{j=m+1}^n n_j \right) a^{(0)}_j =0. $$
Here we use the condition, proved in Section 6, that for sufficiently
small $\epsilon >0$ we have
$$ \la a^{(0)}_j , w_*^0 (a^{(0)}_k) \ra = \delta_{ik}. $$

The above condition can be rewritten as
$$ n_j = - \sum_{i=1}^m n_i N^{W_0}_{k_j}(a^{(1)}_i,a^{(0)}_j)  \ \hbox{ for 
all } \ j=m+1,\ldots, n. $$
Notice that this means that we can choose arbitrarily the first $\{
n_i \}_{i=1,\ldots m}$ coefficients and the remaining $\{ n_j
\}_{j=m+1,\ldots, n}$ are determined by the above relation.

\endproof

We then derive the exactness result as follows.

\begin{Pro}
The inclusion $Ker(w_*^0) \subset Im(w_*^1)$ is satisfied.  
\label{inclusion}
\end{Pro}

\proof
Suppose given an element $a\in Ker(w_*^0)$. By Lemma \ref{kernel}, for 
sufficiently small $\epsilon >0$ in the surgery perturbation, we
can write 
$$ a= \sum_{i=1}^m n_i a^{(1)}_i + \sum_{j=m+1}^n n_j a^{(0)}_j $$
with $n_i$ and $n_j$ in $\Z$, satisfying
$$ n_j = - \sum_{i=1}^m n_i N^{W_0}_{k_j}(a^{(1)}_i,a^{(0)}_j). $$
Consider the element
$$ a_1= \sum_{i=1}^m n_i a_i^{(1)}, $$
with the coefficients $n_i$ as above. We have
$$ w_*^1(a_1)= \sum_{i=1}^m n_i a^{(1)}_i + \sum_{j=m+1}^n \sum_{i=1}^m
n_i N^{W_1}_{\ell_{ij}}(a^{(1)}_i,a^{(0)}_j). $$
The relation
$$ N^{W_1}_{\ell_{ij}}(a^{(1)}_i,a^{(0)}_j) = -
N^{W_0}_{k_j}(a^{(1)}_i,a^{(0)}_j) $$
implies that we have
$$ a= w_*^1(a_1). $$

\endproof

The result of Proposition \ref{inclusion}, together with the relation
$w_*^0\circ w_*^1 =0$ proved in the previous section, implies that we
have exactness in the middle term of the sequence
$$ 0\to C_q (Y_1)\stackrel{w^1_q}{\to} C_q(Y,\mu) \stackrel{w^0_q}{\to}
\bigoplus_{\s_{k}} C_{(q)} (Y_0,\s_{k}) \to 0. $$

\section{The connecting homomorphism}

We have established the exactness of the sequence
$$ 0 \to C_*(Y_1) \stackrel{w_*^1}{\to} C_*(Y,\mu)
\stackrel{w_*^0}{\to}  \oplus_k C_{(*)}(Y_0,\s_k) \to 0. $$
Now we can give a more precise description of the connecting
homomorphism in the induced long exact sequence of Floer
homologies. As we are going to see, the connecting homomorphism
is described in terms of the discrepancy between the boundary operator
$\partial_Y$ of the Floer complex $C_*(Y,\mu)$ and the operator
$\partial_{Y_1}\oplus \partial_{Y_0}$, with $\partial_{Y_0}=\oplus_k
\partial_{Y_0,k}$ on 
$$ C_*(Y_1) \oplus \oplus_k C_{(*)}(Y_0,\s_k). $$

More explicitly, let us identify again the points of $\M_{Y,\mu}^*$ with 
$$ {\cal M}_{Y,\mu}^* = \{ a^{(1)}_i \}_{i=1,\ldots m} \cup \{ a^{(0)}_j
\}_{i=m+1,\ldots n}, $$
as in (\ref{MY}). Consider a cycle $\sum_j x_j a^{(0)}_j$ in $\oplus_k
C_{(*)}(Y_0,\s_k)$. We have 
$$ \partial_{Y_0}(\sum_j x_j a^{(0)}_j)=0. $$
Now consider the corresponding element $\sum_j x_j a^{(0)}_j(\epsilon)$ in
$C_*(Y,\mu)$. We have the following result.

\begin{Lem}
\label{boundary:Y:Y0}
Given a cycle $\sum_j x_j a^{(0)}_j$ in $C_{(*)}(Y_0,\s_k)$, we have
\ba \label{boundary:formula} \begin{array}{l} \partial_{Y,\mu}(\sum_j x_j
a^{(0)}_j(\epsilon))= \\[2mm]
\sum_{i,j} x_j n_{Y,\mu}(a_j^{(0)}(\epsilon),a_i^{(1)}(\epsilon))
a_i^{(1)}(\epsilon) \\[2mm] 
- \sum_{r,j,p} x_j n_{Y,\mu}(a_j^{(0)}(\epsilon),a_r^{(1)}(\epsilon))
N^{W_0}_k(a_r^{(1)}(\epsilon), a_p^{(0)})
a_p^{(0)}(\epsilon). \end{array} \na
Here the first sum is over $a_i^{(1)}$ in $\M_{Y_1}^*$ with
$$ \deg_{Y,\mu}(a_j^{(0)}(\epsilon)) -
\deg_{Y,\mu}(a_i^{(1)}(\epsilon)) =1 $$
and the second sum is over $a_p^{(0)}$ in $\M_{Y_0}(\s_k)$ satisfying 
$$ \deg_{Y_0,\s_k}(a_j^{(0)})-\deg_{Y_0,\s_k}(a_p^{(0)})=1 $$
and $a_r^{(1)}$ in $\M_{Y_1}^*$ satisfying
$$ \deg_{Y,\mu}(a_j^{(0)}(\epsilon)) -
\deg_{Y,\mu}(a_r^{(1)}(\epsilon)) =1 $$
and
$$ \iota(W_0,\s_k,a_i^{(1)}(\epsilon),a_p^{(0)})=0. $$
\end{Lem}

\proof It is sufficient to check that we have
$$ \begin{array}{l} \pi\circ \partial_{Y,\mu}(a_j^{(0)}(\epsilon))=
\partial_{Y_0}(a_j^{(0)}) \\[2mm]
- \sum_{r,j,p} x_j n_{Y,\mu}(a_j^{(0)}(\epsilon),a_r^{(1)}(\epsilon))
N^{W_0}_k(a_r^{(1)}(\epsilon), a_p^{(0)})
a_p^{(0)}(\epsilon), \end{array} $$
where $\pi: C_q(Y,\mu)\to \oplus_k C_{(q)}(Y_0,\s_k)$ is the
projection induced from the identification of the moduli spaces in
\cite{CMW}. 
In order to prove this relation, consider 1-dimensional moduli
spaces $\M^{W_0}_k (a_j^{(0)}(\epsilon), a_p^{(0)})$ and their
compactification. We have boundary strata
$$ \begin{array}{l} \partial\M^{W_0}_k (a_j^{(0)}(\epsilon),
a_p^{(0)})= \\[2mm]
\bigcup_{a_s^{(0)}} \hat
\M_{Y,\mu}(a_j^{(0)}(\epsilon),a_s^{(0)}(\epsilon)) \times \M^{W_0}_k
(a_s^{(0)}(\epsilon), a_p^{(0)}) \cup \\[2mm]
\bigcup_{a_i^{(1)}} \M_{Y,\mu}(a_j^{(0)}(\epsilon),
a_i^{(1)}(\epsilon)) \times \M^{W_0}_k (a_i^{(1)}(\epsilon),a_p^{(0)})
\cup \\[2mm] 
\bigcup_{a_s^{(0)}} \M^{W_0}_k (a_j^{(0)}(\epsilon),a_s^{(0)}) \times
\hat \M_{Y_0,\s_k} (a_s^{(0)},a_p^{(0)}). \end{array} $$
Using the results of Corollary \ref{fund:1} and \ref{fund:2} we can
identify these expressions with
$$ \begin{array}{l}
\M_{Y,\mu}(a_j^{(0)}(\epsilon), a_p^{(0)}(\epsilon)) \cup \\[2mm]
\bigcup_{a_i^{(1)}} \M_{Y,\mu}(a_j^{(0)}(\epsilon),
a_i^{(1)}(\epsilon)) \times \M^{W_0}_k (a_i^{(1)}(\epsilon),a_p^{(0)})
\cup \\[2mm] 
\hat \M_{Y_0,\s_k} (a_j^{(0)},a_p^{(0)}). \end{array} $$
When keeping into account the orientations, this yields the desired
formula. 

\endproof

Using this result we obtain the following.

\begin{Lem}
Suppose given a cycle in $\sum_j x_j a^{(0)}_j$ in
$C_{(*)}(Y_0,\s_k)$. The image of $\sum_j x_j a^{(0)}_j$ under the
connecting homomorphism $\Delta$ is given by
$$ \Delta(\sum_j x_j a^{(0)}_j)= \sum_j x_j
n_{Y,\mu}(a_j^{(0)}(\epsilon), a_i^{(1)}(\epsilon)) a_i^{(1)}, $$
that is, by the algebraic counting of flowlines in 
$$ \hat\M_{Y\times \R,\mu}(a_j^{(0)}(\epsilon),a_i^{(1)}(\epsilon)),
$$
for
$\deg_{Y,\mu}(a_j^{(0)}(\epsilon))-\deg_{Y,\mu}(a_i^{(1)}(\epsilon))=1$.
\label{delta:partial} 
\end{Lem}

\proof Using the result of Theorem \ref{W_1:W_0} we can write
(\ref{boundary:formula}) equivalently as
$$ \begin{array}{l} \partial_{Y,\mu}(\sum_j x_j
a^{(0)}_j(\epsilon)= \\[2mm]
\sum_{i,j} x_j n_{Y,\mu}(a_j^{(0)},a_i^{(1)}) a_i^{(1)}(\epsilon) \\[2mm]
+ \sum_{i,j,p} x_j n_{Y,\mu}(a_j^{(0)}(\epsilon),a_i^{(1)}(\epsilon))
N^{W_1}_\ell(a_i^{(1)}, a_p^{(0)}(\epsilon))
a_p^{(0)}(\epsilon). \end{array} $$
By comparing this expression with
$$ \begin{array}{l} w_*^1( \sum_j x_j n_{Y,\mu}(a_j^{(0)}(\epsilon),
a_i^{(1)}(\epsilon)) a_i^{(1)}) = \\[2mm]
\sum_j x_j n_{Y,\mu}(a_j^{(0)}(\epsilon),
a_i^{(1)}(\epsilon)) a_i^{(1)}(\epsilon) + \\[2mm]
\sum_j x_j n_{Y,\mu}(a_j^{(0)}(\epsilon),
a_i^{(1)}(\epsilon)) N^{W_1}_\ell(a_i^{(1)}, a_p^{(0)}(\epsilon))
a_p^{(0)}(\epsilon), \end{array} $$
we obtain 
$$ \Delta(\sum_j x_j a^{(0)}_j)= 
\sum_j x_j n_{Y,\mu}(a_j^{(0)}(\epsilon), a_i^{(1)}(\epsilon))
a_i^{(1)}. $$
This completes the proof. 

\endproof

\subsection{The surgery triangle}

In this section we give a different description of the 
connecting homomorphism. This will prove 
that, in the $\epsilon \to 0$ limit of the surgery perturbation, the
exact triangle for Seiberg-Witten Floer homology  
is a surgery triangle, that is, the connecting homomorphism in the
exact sequence can also
be described as a map $\bar w^2_*$ induced by a surgery cobordism
$\bar W_2$, and the resulting diagram
$$ C_*(Y_1) \stackrel{w_*^1}{\to} C_*(Y,\mu) \stackrel{w^0_*}{\to}
\oplus_k C_{(*)}(Y_0,\s_k) \stackrel{\bar w^2_*}{\to} C_*(Y_1)[-1] $$
is a distinguished triangle, cf. \cite{Fl}, \cite{BD}, \cite{HYW}.

\begin{Pro}
The connecting homomorphism $\Delta$ in the exact triangle is given by 
the following expression,
$$ \Delta(\sum_j x_j a^{(0)}_j)= \bar w^2_* (\sum_j x_j a^{(0)}_j), $$
for any cycle $\sum_j x_j a^{(0)}_j$ in $C_{(*)}(Y_0,\s_k)$, for some $\s_k$,
where $$\bar w^2_* : \oplus_k C_{(*)}(Y_0,\s_k) \to C_*(Y,\mu)$$
is the homomorphism defined by counting solutions in the
zero-dimensional components of the moduli spaces
$$ \M^{\bar W_2}_k(a_j^{(0)},a_i^{(1)}), $$
over the cobordism $\bar W_2$.
\end{Pro}

\proof We only need to prove that the algebraic counting of flow lines
in $$\hat\M_{Y(r)\times \R,\mu}(a^{(0)}_j(\epsilon),
a_i^{(1)}(\epsilon)),$$
agrees with the algebraic counting of monopoles in the
zero-dimensional components of $\M^{\bar
W_2}_k(a^{(0)}_j,a_i^{(1)})$.

We begin with a few observations on the topology of the
cobordisms. Recall that we have $\bar W = \bar W_0 \cup_Y \bar
W_1$. Moreover, we can write $\bar W= \bar W_2 \# \C P^2$, where the
generator of the homology of the $\C P^2$ summand is the surface 
$S^2=D^2_0 \cup_K D^2_1$, with $D^2_i$ the core handles of the cobordisms
$\bar W_0$ and $\bar W_1$. In our analysis of the geometric limits of
solutions, we have considered the punctured cobordisms $\hat W_1$ and
$\hat W_0$. It is convenient here to consider the punctured
$\overline{ \hat W_1 }=\bar W_1 \backslash \{ x_1 \}$ and $\overline{
\hat W_0 }= \bar W_0 \backslash \{ x_0 \}$, obtained by removing
points $\{ x_i \}$ inside the core disks $D^2_i$. Now consider the
doubly punctured space
$$ \bar W\backslash\{ x_0, x_1 \}\cong \overline{ \hat W_1 }\cup_Y
\overline{\hat W_0 }, $$
with a metric that has cylindrical ends $S^3\times [0,\infty)$ at the
punctures and stretched regions $T^2\times [-r,r]$ inside the standard
Heegaard splittings of these $S^3$ and along the product regions
connecting the solid tori in the Heegaard splittings of $S^3$ to the
tubular neighborhoods of the knot in $Y_1$, $Y$, and $Y_0$.
Then the sphere $S^2=D^2_0 \cup_K D^2_1$ corresponds to a cylinder
$$ S^2 \backslash \{ x_0,x_1 \} \cong S^1 \times \R $$
in $ \bar W\backslash\{ x_0, x_1 \}$.

Thus, when removing the
points $x_0$ and $x_1$, the $\C P^2$ summand in $\bar W$ becomes
a disk bundle ${\cal S}$ over the cylinder $S^1 \times \R$, which, with
respect to a fixed trivialization at the two ends, has Euler number
$+1$. In turn, the doubly punctured cobordism $\bar W_2 \backslash \{
x_0,x_1 \}$ is obtained by replacing this disk bundle with the trivial
disk bundle ${\cal S}_0=D^2 \times S^1 \times \R$. 

\begin{figure}[ht]
\epsfig{file=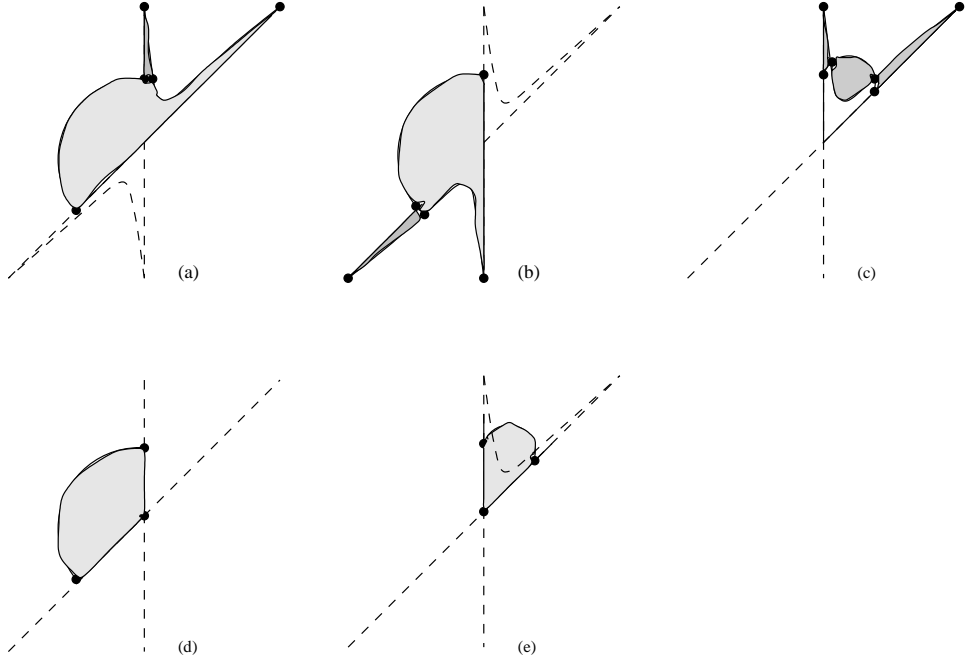,angle=270}
\caption{triangles on $\bar W_2$
\label{figIII:limsW}  }
\end{figure}

Now we proceed to compare the geometric limits of solutions in
$$\M_{Y\times \R,\mu}(a_j^{(0)}(\epsilon), a_i^{(1)}(\epsilon)) \ \
\hbox{ and } \ \ \M^{\bar W_2}_k(a_j^{(0)},a_i^{(1)}).$$ 
There are two distinct cases: 
these correspond to (d) and (e) of Figure \ref{figIII:limsW}.

\noindent {\bf Case 1.}
Consider geometric limits of solutions on the doubly punctured
cobordism $\bar W\backslash\{ x_0, x_1 \}$. By the previous analysis
of the geometric limits of solutions on the punctured cobordisms $W_1$
and $W_0$, together with the gluing theorem for 
$ \M^{\bar W}_{\ell,k}(a_j^{(0)}, a_i^{(1)})$, which in this case gives
$$ \begin{array}{l} \M^{\bar W}_{\ell,k}(a_j^{(0)}, a_i^{(1)})\cong \\[2mm]
\M^{\bar W_0}_k(a_j^{(0)},a_i^{(1)}(\epsilon))\times \M^{\bar
W_1}_\ell(a_i^{(1)}(\epsilon),a_i^{(1)}) \cup \\[2mm]
\M^{\bar W_0}_k(a_j^{(0)},a_j^{(0)} (\epsilon))\times \M^{\bar
W_1}_\ell(a_j^{(0)} (\epsilon),a_i^{(1)}) \end{array} $$ 
for zero-dimensional moduli spaces, we obtain
geometric limits as in 
Figure \ref{figIII:limsW}, (a), (b). These reduce to solutions on the
once punctured $\bar W_2\backslash \{ x_2 \}$ that have geometric
limits as in Figure \ref{figIII:limsW} (d). The counting
$$ N^{\bar W_0}_k(a_j^{(0)},a_i^{(1)}(\epsilon)) + N^{\bar
W_1}_\ell(a_j^{(0)} (\epsilon),a_i^{(1)})=0 $$
implies that the total counting of these solutions in the
zero-dimensional moduli space $\M^{\bar W_2}_k (a_j^{(0)},a_i^{(1)})$
is zero. Notice that, in this case we have 
$$ \M_{Y(r)\times \R,\mu}(a_j^{(0)}(\epsilon), a_i^{(1)}(\epsilon)) =
\emptyset, $$
therefore we also have $n_{Y,\mu}(a_j^{(0)}(\epsilon),
a_i^{(1)}(\epsilon))=0$. 

\noindent {\bf Case 2.}
In this case,  the  
zero-dimensional moduli space $\M^{\bar W_2}_k (a_j^{(0)},a_i^{(1)})$
is obtained from approximate solutions of the form
$$ \hat \M_{V\times \R}([A,\psi],[\tilde A,\tilde \psi], a_\infty) \#
\Delta, $$
with $\Delta$ a triangle as in Figure \ref{figIII:limsW} (e), and 
$([A,\psi],[\tilde A,\tilde \psi], a_\infty)$ is the
admissible triple with respect to $(a_j^{(0)},a_i^{(1)})$. 
(The corresponding holomorphic triangle for geometric limits
on $\bar W -\{ x_0, x_1\}$ is illustrated in Figure 4 (c).) 
This means
that, in forming the preglued solutions, elements in the moduli spaces
$$ \M_{V\times \R}([A,\psi],[\tilde A,\tilde \psi], a_\infty) $$
that differ by a translation will give rise to the same solution on
$\bar W_2$. Notice that geometric limits for 
 $\M_{Y(r)\times
\R}(a_j^{(0)}(\epsilon),a_i^{(1)}(\epsilon))$ are given by
$$ \M_{V\times \R}([A,\psi]^\epsilon,[\tilde A,\tilde \psi]^\epsilon, a_\infty) \#
\Delta^\epsilon $$
with $\Delta^\epsilon$ a holomorphic disk in $H^1(T^2,\R)$ with boundary along
arcs of $\ell_\mu$ and $\ell$ connecting the asymptotic limits
of $a_j^{(0)}(\epsilon),a_i^{(1)}(\epsilon)$ respectively,  the 
corresponding admissible triple is 
$([A,\psi]^\epsilon,[\tilde A,\tilde \psi]^\epsilon, a_\infty)$
for $(a_j^{(0)}(\epsilon),a_i^{(1)}(\epsilon))$.

As we let $\epsilon \to 0$ these geometric limits approximate
geometric limits
$$ \hat \M_{V\times \R}([A,\psi],[\tilde A,\tilde \psi], a_\infty) \#
\Delta$$
of  $\M^{\bar W_2}_k (a_j^{(0)},a_i^{(1)})$.
Under the gluing map, we obtain that the
algebraic counting of solutions in the zero-dimensional moduli space
$\M^{\bar W_2}_k (a_j^{(0)},a_i^{(1)})$ agrees with 
 the algebraic counting of
flow lines in $\hat\M_{Y\times
\R,\mu}(a_j^{(0)}(\epsilon),a_i^{(1)}(\epsilon))$, 
$$ N^{\bar W_2}_k(a_j^{(0)},a_i^{(1)})=
n_{Y,\mu}(a_j^{(0)}(\epsilon),a_i^{(1)}(\epsilon)) . $$
This completes the proof.

\endproof

Let us recall briefly the following definitions of homological
algebra, cf. \cite{GM}.

A triangle (in a category of complexes) is a diagram of the form
$$ K_\bullet \stackrel{u}{\to} L_\bullet \stackrel{v}{\to} M_\bullet
\stackrel{w}{\to} K_\bullet [-1]. $$
Here we take complexes with differentials of degree $(-1)$ for
consistency with our notation. $K_\bullet [n]$ is the $n$-th shift
with $K_i[n]=K_{i+n}$, and $d_{K[n]}=(-1)^n d_K$. A triangle is {\em
distinguished} if it is {\em quasi--isomorphic} to a diagram of the
form 
$$ K_\bullet \stackrel{\tilde u}{\to} \hbox{Cyl}(u)\stackrel{\pi}{\to}
\hbox{C}(u) \stackrel{\delta}{\to} K_\bullet [-1]. $$
Here $\hbox{Cyl}(u)$ and $\hbox{C}(u)$ are the cylinder and the cone
on the morphism $u$, defined as in \cite{GM} pg.154 (up to adjusting
notations to differentials of degree $(-1)$), and $\delta$ is the
connecting homomorphism. 

With this notation and the previous results, we can rephrase our main
result as the following Proposition. 

\begin{Pro}
For sufficiently small $\epsilon \leq \epsilon_0$ in the surgery
perturbation $\mu$, the sequence
$$ C_*(Y_1) \stackrel{w_*^1}{\to} C_*(Y,\mu) \stackrel{w^0_*}{\to}
\oplus_k C_{(*)}(Y_0,\s_k) \stackrel{\bar w^2_*}{\to} C_*(Y_1)[-1] $$
is a distinguished triangle.
\label{surgery-triangle}
\end{Pro}

\proof
The result follows from the previous results, together with
Proposition 5, pg.157 of \cite{GM}. In fact, we have proven that, as
$\epsilon \to 0$, the sequence 
$$ 0 \to C_*(Y_1) \stackrel{w_*^1}{\to} C_*(Y,\mu)
\stackrel{w^0_*}{\to}\oplus_k C_{(*)}(Y_0,\s_k) \to 0 $$
is an exact triple of complexes (\cite{GM} pg. 42). Proposition 5,
pg.157 of \cite{GM} then gives the result.

\endproof

\noindent {\bf Matilde Marcolli}, 
Max--Planck--Institut f\"ur Mathematik, D-53111 Bonn, Germany. 
marcolli\@@mpim-bonn.mpg.de

\vskip .2in

\noindent {\bf Bai-Ling Wang}, Department of Pure
Mathematics, University of Adelaide, Adelaide SA 5005. 
 bwang\@@maths.adelaide.edu.au \par
\noindent Max--Planck--Institut f\"ur Mathematik, D-53111 Bonn, Germany.
\par \noindent bwang\@@mpim-bonn.mpg.de

\end{document}